\documentclass[12pt]{article}

\usepackage[utf8]{inputenc}
\usepackage[top=1in,bottom =1in, left=1in,right=1in]{geometry}
\usepackage{amsmath,amssymb,mathtools,bm,isomath,amsthm}
\usepackage{graphicx}
\usepackage{epstopdf}
\usepackage{color}
\usepackage[dvipsnames]{xcolor}
\usepackage{soul}
\usepackage{enumerate}
\usepackage{nicefrac}
\usepackage{enumitem}
\usepackage{mathrsfs}
\usepackage{scalerel,stackengine}
\usepackage{afterpage}
\usepackage{graphicx}
\usepackage{url}
\usepackage[affil-it]{authblk} 
\usepackage[toc,page]{appendix}
\usepackage{hyperref}
\usepackage{amsmath}
\usepackage{tikz}
\usepackage{float}
\usepackage{optidef}
\usepackage{algpseudocode}
\usetikzlibrary{positioning}
\usetikzlibrary{positioning, calc, decorations.pathreplacing}
\usetikzlibrary{matrix, positioning}
\usetikzlibrary{fit,positioning}
\usepackage[export]{adjustbox}
\graphicspath{{./images/}{./dia/}{./}}
\usepackage{pgfplots}
\usepackage{graphicx}
\usepackage{epstopdf}  
\usepackage{flafter,placeins}
\usepackage{subcaption} % figures & sous figures / figures & subfigures
\usepackage{listofitems}
\usepackage{tikz}
\usepackage{siunitx} % unites SI / SI units
\usepackage{amssymb} % a
\stackMath
\newcommand\reallywidehat[1]{%
	\savestack{\tmpbox}{\stretchto{%
			\scaleto{%
				\scalerel*[\widthof{\ensuremath{#1}}]{\kern-.6pt\bigwedge\kern-.6pt}%
				{\rule[-\textheight/2]{1ex}{\textheight}}%WIDTH-LIMITED BIG WEDGE
			}{\textheight}% 
		}{0.5ex}}%
	\stackon[1pt]{#1}{\tmpbox}%
}

\newcommand{\x}{\boldsymbol{x}}

\pgfplotsset{compat=1.18} 
\begin{document}
	\setlength{\baselineskip}{15pt}
	
	\title{Coupling of the Finite Element Method with Physics Informed Neural Networks for the Multi-Fluid Flow Problem }
	
	\author[1]{Michel Nohra}
	\author[1]{Steven Dufour}
	\affil[1]{D\'epartement de math\'ematiques et de g\'enie industriel, Polytechnique Montr\'eal, Montr\'eal, Qu\'ebec, Canada, H3T 1J4}
	
	\date{\today}
	
	\maketitle

\begin{abstract}

Multi-fluid flows are found in various industrial processes, including metal injection molding and 3D printing. The accuracy of multi-fluid flow modeling is determined by how well interfaces and capillary forces are represented. In this paper, the multi-fluid flow problem is discretized using a combination of a Physics-Informed Neural Network (PINN) with a finite element discretization. To determine the best PINN formulation, a comparative study is conducted using a manufactured solution. We compare interface reinitialization methods to determine the most suitable approach for our discretization strategy. We devise a neural network architecture that better handles complex free surface topologies. Finally, the coupled numerical strategy is used to model a rising bubble problem.
\end{abstract}

  \section{Introduction}

Multi-fluid flows play an important role in various industrial processes, ranging from 3D printing to oil extraction processes, and metal injection applications. Understanding and accurately modeling these complex flows is important for optimizing process efficiency, and for ensuring reliable outcomes. The numerical modeling of multi-fluid flows is challenging due to the dynamics of the interfaces between the different fluids, and the influence of capillary force.
Numerical methods for the study of multi-fluid flows can be broadly classified into Eulerian and Lagrangian techniques. The Lagrangian approach splits the domain into subdomains and discretizes the interfaces. The nodes on each interface are transported by the fluid flow. This process necessitates updating physical quantities at every time step, and it requires the remeshing of the entire domain if elements become too degenerated. Geometric quantities, such as the curvature of the interface and the normal vectors to the interfaces, are obtained from the relative position of these nodes, which may become problematic for a distorted discretization of the interface. The capillary force is treated as a boundary condition for each subdomain.

On the other hand, Eulerian techniques employ a different strategy for modeling multi-fluid flows. A marker variable is used to locate each fluid, and the evolution of the interfaces is governed by the mass conservation law
$$\partial_t F + \boldsymbol{u} \cdot \nabla F =0, $$
where $\boldsymbol{u}$ denotes the velocity of each fluid, and $F$ is the marker. Eulerian methods treat fluid properties as continuous fields, defined in the entire computational domain. Geometric properties, such as the curvature and normal vectors to an interface, are calculated using the marker, where normal vectors become $\boldsymbol{n}=\frac{\nabla F}{\| \nabla F \|}$ and the curvature of the interface is given by $\kappa=\nabla \cdot \boldsymbol{n} $, thus approximating the capillary forces by a volumetric force. Eulerian techniques are more robust and more suitable for distorted $3$D interfaces, but suffer from numerical errors related to the dissipation of the interface marker. 
With the Eulerian approach, three primary methods are commonly used to simulate multi-fluid flows: the volume-of-fluid (VOF) method, the level-set method, and the pseudo-concentration method. Each method offers unique features and advantages to accurately represent and capture fluid interfaces.

With the volume-of-fluid method, the marker takes values between zero and one, representing the fraction of a given fluid within a computational cell. If the marker value in a cell is $F=1$, then the cell contains fluid $\text{A}$; if the value is $F=0$, the cell contains only fluid $\text{B}$. This marker is then transported using a finite volume method. The cells containing an interface are such that $0<F<1$, but the position of the interface inside the cell is unknown. After transporting the marker, algorithms are used to reconstruct an interface and extract geometric quantities such as the normal vectors and the curvature of an interface. 

The level-set method, on the other hand, represents the location of the two fluids using a signed distance function to an interface. 
%This results in a smooth marker with a unitary gradient. 
The zero contour of the marker determines the location of an interface. The normal vectors are given by $\boldsymbol{n}=\frac{\nabla F}{\| \nabla F \|},$ and the curvature of interfaces by $\kappa=\nabla \cdot \boldsymbol{n}.$ 
The pseudo-concentration method, sometimes referred to as the conservative level-set method, combines components of both the level-set and the volume-of-fluid methods. This method uses a distance function to an interface to approximate a smooth Heaviside function that varies between zero and one, which is used as a marker. The marker's contour level of $0.5$ determines each interface position, and the geometric quantities are obtained using the expressions used with the level-set method.
Notably, the pseudo-concentration method exhibits a steeper gradient at the interface while maintaining a constant value away from it. 
In this work, we will use the level-set method due to the smoothness of the marker.

Traditional numerical methods, like the finite difference, the finite element, and the finite volume methods, are commonly used to numerically model multi-fluid flow problems. The discretizations employed by these methods result in high computational and memory requirements for high-dimensional problems. In addition, the second-order derivatives required to calculate the curvature are a source of additional numerical errors.

Neural networks are revolutionary mathematical approximators that have been effectively applied in various fields, including speech recognition and classification tasks \cite{nn_1,nn_2,nn_3}.
With recent advancements in computer hardware, the use of neural networks in scientific computing has recently grown in popularity. These neural networks offer promising alternatives for numerical modeling and simulation tasks, leveraging their ability to learn complex patterns and relationships from large datasets.
Neural networks have been applied to multiple scientific problems, such as turbulence model design from data \cite{turb1, turb2,turb3,turb4} and for simulating heat transfer problems \cite{HT1,ht2}.

Among the methods that use deep learning for scientific computing, Physics-Informed Neural Networks (PINNs)\cite{raissi} are a popular approach for simulating physical problems. PINNs, first introduced by Raissi et al.\cite{raissi}, leverage all the information available about the problem, such as the underlying PDE, boundary conditions, physical properties, and known measurement data, to train a neural network that best approximates the model. These networks transform boundary-value problems into optimization problems, providing flexibility to include any constraint in the problem.

Since their introduction, PINNs have been applied to a variety of forward, inverse, and parameterized problems \cite{HT1,nsfnets,ht2,max_net,max2,Max1,GAPINN}, including fluid dynamics problems. However, due to the second-order, non-linear and saddle point nature of fluid dynamics problems, traditional numerical methods still outperform PINNs for forward problems (well-posed boundary-value problems). PINNs become interesting when experimental data are available \cite{PINN_fluid}.

In the case of hyperbolic equations, such as the transport equation, traditional methods, such as the finite element method, often exhibit instabilities that necessitate special treatment. This is where PINNs offer an advantage, as they are not prone to the same instabilities, and they do not require additional measures to ensure numerical stability. Therefore, PINNs can be seen as a valuable alternative for the accurate discretization of hyperbolic equations. Another advantage of PINNs over traditional methods is that the numerical approximations are continuous and infinitely differentiable, which is important for an accurate approximation of surface tension.
For the reasons listed above, Physics-Informed Neural Networks will be used for the discretization of the transport equation of the level-set function, and the finite element method will be used to discretize the Navier-Stokes equation.

In this work, we will describe the multi-fluid flow problem and the governing equations, the FEM for the Navier-Stokes equations, and the PINNs. Additionally, we will propose three different time-integration schemes for the advection equation with PINNs. We will
present a new approach, based on the multi-level PINN~\cite{ziad}, to obtain a more accurate approximation of the level-set function. We introduce three different methods to reinitialize the level-set function. And finally, validate the proposed methodologies on various problems with manufactured solutions, and on the rising bubble problems.

\section{Governing equations}
Let us consider a multi-fluid flow problem with two non-miscible and incompressible fluids $\text{A}$ and $\text{B}$. The density and viscosity of each fluid are given by $\rho_\text{A}, \mu_\text{A}$, $\rho_\text{B}$ and $\mu_\text{B}$, and are defined over subdomains $\upOmega_\text{A}$ and $\upOmega_\text{B}$, where $\upOmega_\text{A} \cup \upOmega_\text{B} = \upOmega$, which is the whole domain, and  $\upOmega_\text{A} \cap \upOmega_\text{B} = {\mathrm{\Gamma}},$ where $\mathrm{\Gamma}$ is the interface between the two fluids.
With the level-set method, we define a signed distance function to the interface $\mathrm{\Gamma},$ denoted by $F$, where $F$ is negative in $\upOmega_\text{A},$ positive in $\upOmega_\text{B}$, and zero at the interface. We can define an approximate Heaviside function using the level-set function $F$. Various approximations can be used, including the sigmoid function $H(F)=\sigma(\alpha F)$, where $\alpha$ determines the sharpness of $H$, which depends on the grid size. Using the level-set function, and the approximate Heaviside function $H$, we can build a one-fluid model with variable physical quantities:
\begin{equation*}
    \begin{aligned}
        \rho(F) = \rho_\text{A} H(F) + \rho_\text{B} (1-H(F));\\
        \mu(F) = \mu_\text{A} H(F) + \mu_\text{B} (1-H(F)).
    \end{aligned}
\end{equation*}
%$$P = P_\text{A} H + P_\text{B} (1-H), $$
%where $P$ defined on all $\upOmega$ represent any physical quantity.
Using this one-fluid model, we obtain the set of equations that represent the continuity or incompressibility equation, the conservation of momentum equations, and the marker advection equation:
\begin{equation}
\begin{aligned}
& \quad \nabla \cdot \boldsymbol{u} = 0; \\
& \quad \rho(F) \frac{\partial \boldsymbol{u}}{\partial t} + \rho(F) (\boldsymbol{u} \cdot \nabla)\boldsymbol{u}  + \nabla p -  \nabla \cdot \nabla( \mu(F)\boldsymbol{u}) = \boldsymbol{F}_\text{ext};\\
& \quad \partial_t F + \boldsymbol{u} \cdot \nabla F =0,
\end{aligned}
\end{equation}
 where $\boldsymbol{u}$ represents the fluids velocity, $p$ represents their pressure, $F$ is the level-set function, $\boldsymbol{F}_{\text{ext}}$ represents external forces, $\mu$ is the viscosity and $\rho$ is the density of the fluids.
 
$\boldsymbol{F}_{\text{ext}}$ includes an approximation of the capillary force, which is a surface force at the interface,
$$\boldsymbol{f_\text{s}} = \delta \gamma \kappa \boldsymbol{n}_\mathrm{\Gamma},$$
where $\delta$ is the Dirac function defined at the interface, $\gamma$ is the surface tension coefficient, $\kappa$ is the interface curvature and $\boldsymbol{n}_\mathrm{\Gamma}$ is a normal vector at the interface. Using the continuum surface force approximation \cite{brackbill1992continuum}, we can approximate the normal vectors using $\boldsymbol{n}_\mathrm{\Gamma}\ = \frac{\nabla F}{\|\nabla F\|}$, the Dirac function by $\delta=\frac{\partial H}{\partial F}$, and the curvature of the interface by $\kappa = \nabla \cdot \boldsymbol{n}_\mathrm{\Gamma},$ and obtain a volumetric approximation of the capillary force.

\subsection{Adimensional equations}
Adimensionalizing, often referred to as non-dimensionalization, is a fundamental technique in the fields of physics and engineering that helps reduce the number of parameters in a problem, allowing for a clearer understanding of the underlying dynamics and relationships. Results derived from dimensionless equations are often more general, as they are not tied to specific units of measurement. This universality facilitates easier comparison across different systems and scales, leading to broader insights and more robust solutions.  Furthermore, adimensionalizing the problem means that the calculated variables would tend to be in the order of unity, making the problem easier from a numerical point of view.
The quantities of interest would be reformulated as follows:
\begin{table}[ht!]
    \centering
    \begin{tabular}{|c|c|c|c|} \hline 
         $\rho = \hat{\rho} \bar{\rho} $&  $\mu = \hat{\mu} \bar{\mu} $&  $\boldsymbol{u} = \hat{u} \bar{\boldsymbol{u}} $&  $L = \hat{L} \bar{L} $\\ \hline 
         $t=\frac{\hat{L}}{\hat{u}}\bar{t}$&  $\partial_t = \frac{\hat{u}}{\hat{L}} \partial_{\bar{t}}$&  $\nabla = 1/ \hat{L} \bar{\nabla} $&  $p = \hat{\rho} \, \hat{u}^2 \, \bar{p}$\\ \hline
    \end{tabular}
    \caption{Adimensional properties.}
\end{table}
 where $\bar{A}$ represent the non-dimensional value of $A$ and $\hat{A}$ the characteristic value of $A$.
 In the remainder of this work, all variables will be presented in non-dimensional form, and the hat notation ($\hat{\cdot}$) will be omitted for simplicity.
After substituting in the differential equation, we obtain the following non-dimensional numbers:
\begin{itemize}
    \item Reynolds number $Re = \frac{\rho u L}{\mu}$ and it signifies the importance of inertial forces over viscous forces. A higher Reynolds number means the flow exhibits more turbulent behavior, and from a computational point of view, it means that our problem has less diffusion, and necessitates a finer mesh for a stable solution.
    \item Weber number $We = \frac{\rho u^2 L}{\sigma}$ and it represents the importance of inertial forces over capillary forces. A higher Weber number means that the interface has a greater tendency to break and contains higher curvatures. 
    \item Froud number $Fr= \frac{u}{\sqrt{gL}}$ and represent the importance of inertial forces over gravitational forces.
\end{itemize}
The complete non-dimensional equations become: 
    
\begin{equation}
\begin{aligned}
&\nabla \cdot \boldsymbol{u} =0;
\\
&\rho(F) \frac{\partial \boldsymbol{u}}{\partial t} + \rho(F) (\boldsymbol{u} \cdot \nabla)\boldsymbol{u}  + \nabla p - \frac{1}{Re} \mu(F) \nabla^2 \boldsymbol{u} =\frac{1}{We} \frac{\nabla F}{|| \nabla F ||} \nabla \cdot (\frac{\nabla F}{|| \nabla F ||} ) + \frac{1}{Fr^2}\rho(F) ;
\\
&\partial_t F + \boldsymbol{u} \cdot \nabla F =0.
\end{aligned}
\end{equation}

\section{Numerical methods}
In this work, we combine two numerical methods: the Finite Element Method (FEM) for the discretization of the fluid dynamics problem, to determine the fluid's velocity and pressure, and the Physics-Informed Neural Network (PINN) for the advection of the level-set function.

 \subsection{The finite element method for the Navier-Stokes equations}

 The level-set method for modeling multi-fluid flows involves the numerical discretization of nonlinear, coupled partial differential equations. This process typically progresses through discrete time steps. At each time step, the Navier-Stokes equations are discretized based on the previously computed level-set function. The level-set function is then advected using the computed velocity field.

In this section, we explain the discretization of the Navier-Stokes equations using the finite element method. Let us define the functional space of square-integrable functions
$$L^2(\upOmega) = \{ u : \upOmega \xrightarrow{}\mathbb{R} | \int_\upOmega u^2 dx < \infty \},$$ 
with the scalar product 
$$(u,v)_{0,\upOmega} = \int_\upOmega u \, v \, dx.$$
The vectorial and tensorial versions are 
$$L^2(\upOmega)^n = \{ \boldsymbol{u} : \upOmega \xrightarrow{}\mathbb{R}^n | \int_\upOmega \boldsymbol{u} \cdot \boldsymbol{u} \ dx < \infty \},$$ 
and
$$L^2(\upOmega)^{(n \times n)} = \{ \boldsymbol{\psi} : \upOmega \xrightarrow{}\mathbb{R}^{(n \times n)} | \int_\upOmega \boldsymbol{\psi} : \boldsymbol{\psi} \ dx < \infty \},$$
where '$:$' is the double contraction operator.
Their respective scalar product are $$(\boldsymbol{u},\boldsymbol{v})_{0,\upOmega} = \int_\upOmega \boldsymbol{u} \cdot \boldsymbol{v} \ dx,$$ and 
$$(\boldsymbol{\psi},\boldsymbol{\phi})_{0,\upOmega} = \int_\upOmega \boldsymbol{\psi} : \boldsymbol{\phi} \ dx.$$
Next, we define the Sobolev spaces $$ H^1(\upOmega) = \{ u \in L^2(\upOmega) \ | \ \partial_{x_i} u \in L^2(\upOmega), i=1,2,...,n \} ,  $$
and 
$$ H^1(\upOmega)^n = \{ \boldsymbol{u} \in L^2(\upOmega)^n \ | \ \partial_{x_i} \boldsymbol{u}_j \in L^2(\upOmega), \ i,j=1,2,...,n \} ,  $$
with the inner product  
$$(u,v)_{1,\upOmega} = (u , v)_{L^2(\upOmega)} + \sum_{i=1}^n (\partial_{x_i}u , \partial_{x_i}v )_{L^2(\upOmega)}, $$
and finally, for imposing boundary conditions, we define the Sobolev space
\[H^1_{\text{D} \boldsymbol{u}}(\upOmega)^n = \{\boldsymbol{u} \in H^1(\upOmega)^n \ |\ \boldsymbol{u}(\x)=0 \ \text{for} \ \x \in \mathrm{\Gamma}_{\text{D} \boldsymbol{u}} \} ,\]

and \[ H^{1/2}(\mathrm{\Gamma}_0)^n = \{ \gamma_0 (\boldsymbol{u}) \in L^2({\rm \Gamma}_0)^n |\ \boldsymbol{u} \in H^1(\upOmega)^n \}  ,\]
where $\gamma_0(\boldsymbol{u})$ represents the trace of $\boldsymbol{u}$ on ${\rm \Gamma}_0$.

The weak formulation of the equations is obtained by multiplying these equations by test functions belonging to appropriate functional spaces, and integrating by parts over $\upOmega$. The mass conservation equation becomes
\begin{equation}
( q , \nabla \cdot \boldsymbol{u} )_{0,\upOmega}= 0 ,\ \ \ \forall q \in L^2(\upOmega).
\label{eqn:weak_mass}
\end{equation}
Given $\boldsymbol{f} \in L^2(\upOmega)^n $, $\boldsymbol{u}_{\text{D} \boldsymbol{u}} \in H^{1/2}(\mathrm{\Gamma}_{\text{D} \boldsymbol{u}})^n $ and $\boldsymbol{t}_{\text{N} \boldsymbol{u}} \in H^{-1/2}(\mathrm{\Gamma}_{\text{N} \boldsymbol{u} })$, where $H^{-1/2}$ is the dual space of $H^{1/2}$, we obtain the weak formulation of the momentum equation, 
\begin{align}\label{eqn:weak_NSt}
 \big(\rho \frac{\partial \boldsymbol{u}}{\partial t} , \boldsymbol{v}\big)_{0,\upOmega} + (\mu \nabla \boldsymbol{u} , \nabla \boldsymbol{v})_{0,\upOmega} 
 +& (\rho (\boldsymbol{u} \cdot \nabla) \boldsymbol{u} , \boldsymbol{v})_{0,\upOmega} - ( p , \nabla \cdot \boldsymbol{v})_{0,\upOmega}\\
 =&
\ (\boldsymbol{f} , \boldsymbol{v})_{0,\upOmega} + \langle t_{\text{N} \boldsymbol{u}} , v\rangle_{H^{-1/2} ,H^{1/2} },  \quad  \forall \: \boldsymbol{v} \in H^1_{\text{D} \boldsymbol{u}}(\upOmega)^n \notag  .
\end{align}

The bilinear form obtained for the weak formulations in equations \eqref{eqn:weak_mass} and \eqref{eqn:weak_NSt} can be written as:
\begin{align*}
A(\boldsymbol{u},\boldsymbol{v}) +B (\boldsymbol{v},p) &= G(\boldsymbol{v}) \ , \ \  \forall  \boldsymbol{v} \in V ;\\
B(\boldsymbol{u},q) &= 0 \ , \ \  \forall  q \in Q,
\end{align*}
where $V$ and $Q$ are Hilbert spaces, and
\(A\) and \(B\) are bilinear forms given by
\begin{equation*}
    \begin{aligned}
&A(\boldsymbol{u},\boldsymbol{v}) = \big(\rho \frac{\partial \boldsymbol{u}}{\partial t} , \boldsymbol{v}\big)_{0,\upOmega} + (\mu \nabla \boldsymbol{u} , \nabla \boldsymbol{v})_{0,\upOmega} 
 + (\rho (\boldsymbol{u} \cdot \nabla) \boldsymbol{u} , \boldsymbol{v})_{0,\upOmega} ;   \\
&B (\boldsymbol{v},p) = - ( p , \nabla \cdot \boldsymbol{v})_{0,\upOmega},
    \end{aligned}
\end{equation*}
 
and \(G\) is a linear form, \[G(v) = (\boldsymbol{f} , \boldsymbol{v})_{0,\upOmega}. \] 
A unique solution exists according to the Lax-Milgram theorem, and the Ladyzhenskaya-Babuska-Brezzi (LBB) condition is satisfied, if the following conditions are met:
\begin{itemize}
\item Coercivity: This refers to the property that the bilinear form $A(u, v)$ is coercive on the solution space V:
\[ \exists \alpha > 0 \text{ such that } \forall\: u \in V , \ A(u, u) \geq \alpha \|u\|_V^2;
\]
\item Continuity (or boundedness): This condition ensures that the bilinear form is continuous:
$$|A(u, v)| \leq M \|u\|_V \|v\|_V;$$
$$|B(u, q)| \leq N \|u\|_V \|q\|_Q,$$
for all $u, v \in V$,   $q \in Q$, for positive constants $M$ and $N$.
\item LBB (inf-sup) condition:  Ensuring that the chosen spaces for velocity and pressure are compatible is crucial. This compatibility condition is often expressed via the inf-sup condition,
$$ \inf_{v \in V} \sup_{q \in Q} \frac{ |B(v, q)|}{\|v\|_V \|q\|_Q} > \beta,$$
where $\beta$ is a positive constant.
\end{itemize}

In the finite element setting, the domain $\upOmega$ is discretized with smaller subdomains, called elements.
We take the discrete functional space $V_h \subset H^1_{\text{D} \boldsymbol{u}}(\upOmega)^n $ and $Q_h \subset L^2(\upOmega)$ of finite dimensions, with basis functions $\boldsymbol{\phi}_i$ and $\eta_i$ respectively, defined on each element, to express a discrete approximation of the dependent variables as:
\begin{align*}
\boldsymbol{u}_h &= \sum_{i=1}^{N_u}u_i\boldsymbol{\phi}_i(x,y); \\
p_h &= \sum_{i=1}^{N_p}p_i\eta_i(x,y),
\end{align*}
where $\boldsymbol{u}_h$ is the discrete approximation of the velocity field, and $p_h$ is the discrete approximation of the pressure. With the Galerkin finite element method, we consider $\boldsymbol{v}_h \in V_h $ and $q_h \in Q_h$.
As for the choice of the subspaces $Q_h$ and $V_h$, they need to satisfy the LBB condition defined above. An example of such spaces is the Taylor-Hood element combination, where the \(\boldsymbol{\phi}_i\) are the quadratic Lagrange elements \(P_2\), and the \(\eta_i\) are the linear Lagrange elements \(P_{1}\).

\subsection{Physics-Informed Neural Networks}

% To gain a better understanding of neural networks, let us begin with a simple example of a one-layer, single-neuron neural network. In this case, we denote the input and output as $x$ and $x_1$, respectively. The equation governing the output $x_1$ can be written as 
% $$x_1 = \sigma (W x + b),$$
% where $W$ and $b$ are the weight and bias of the neuron. The activation function $\sigma$ is chosen by the user and can take various forms, such as cosine or hyperbolic tangent.
To gain a better understanding of neural networks, let us begin with a simple example of a one-layer neural network. If we consider the input vector $\x$ and the output vector $\x_1$ to have dimensions $n$ and $m$ respectively, the output can be expressed as $$\boldsymbol{x}_1 = \sigma (W \boldsymbol{x} + \boldsymbol{b}),$$ where $W$ is a weight matrix of dimension $m \times n$ and $\boldsymbol{b}$ is a bias vector of dimension $m$. The activation function $\sigma(\boldsymbol{z})$  is user-defined and can take various forms, such as a cosine or a hyperbolic tangent. This function is applied to each component of a vector $\boldsymbol{z}$.  

\begin{figure}
    \centering
    \includegraphics{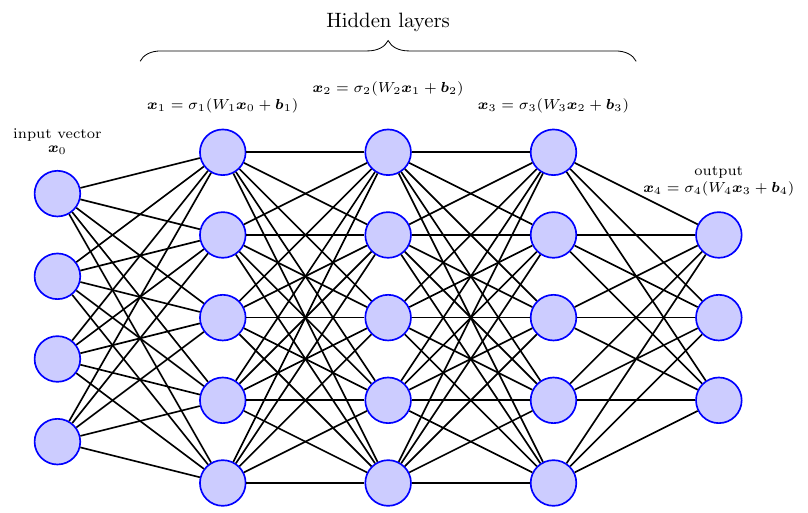}
    \caption{Feedforward neural network}
    \label{fig:FFNN_T1}
    
\end{figure}

 By stacking multiple layers, we can establish a feedforward neural network, which involves passing the output of one layer to the input of the next layer, as illustrated in figure~\ref{fig:FFNN_T1}. The output of the network is then expressed as 
\[ N(\x) = \sigma_l \big( W_n  \sigma_{l-1}( W_{l-1} \sigma_{l-2}(\ldots (\sigma_1(W_1 \x + \boldsymbol{b}_1  )  ) \ldots  ) +\boldsymbol{b}_{l-1}  )       +\boldsymbol{b}_l     \big).  \]

One can build a more complex neural network by changing the connections (the architecture) between the layers.

Once the neural network has been constructed with the desired architecture, the next step is to determine the optimal values for its parameters $\boldsymbol{\zeta},$ a vector that contains all the weights and the biases of the layers of the neural network, so that it produces the desired output. This process of fine-tuning the weights is referred to as ``training''. Training involves minimizing the disparity between the actual output generated by the network and the desired output. To formulate this as an optimization problem, we aim to minimize a specific loss function, denoted as $L$, which is user-defined. For example, for regression problems, the Mean Squared Error (MSE) is commonly used as a loss function, while for classification tasks, the Cross-Entropy loss function is often preferred.
The optimization problem can then be expressed as
\begin{equation}    	\min_{\boldsymbol{\zeta} }\  L(N_{\boldsymbol{\zeta}}(\x),D(\x)),
\end{equation}
where $N_{\boldsymbol{\zeta}}(\x)$ represents the output produced by the neural network for a given ${\boldsymbol{\zeta}}$.
$D(\x)$ is the target, or desired output.
$L$ is the user-defined loss function that quantifies the dissimilarity between the predicted and the actual outputs.
While various optimization methods could address this minimization problem, in the context of this work, we will only consider gradient-based optimization methods.
Gradient-based methods leverage the property that the gradient of a function points in the direction of the function's steepest ascent. Therefore, to minimize the function, one would generally move in the direction opposite to the gradient, with a step size proportional to a user-defined learning rate.

Popular gradient-based methods include the Stochastic Gradient Descent (SGD), the Limited-memory Broyden-Fletcher-Goldfarb-Shanno method (L-BFGS) \cite{bollapragada2018progressive}, and the Adaptive Moment Estimation method (Adam) \cite{adam}. With neural networks, the gradient is calculated using the back-propagation algorithm, which employs the chain rule and automatic differentiation to compute the partial derivatives of the network's output with respect to the weights, the biases, and the network's input.

With the original PINN, introduced by Raissi et al. \cite{raissi}, which is usually referred to as the vanilla PINN, we approximate the solution of a boundary-value problem using a feedforward neural network. The derivatives of the approximation are obtained by automatic differentiation. 
Let us consider a domain $\upOmega$ with a boundary $\mathrm{\Gamma} = \partial \upOmega $. We want to find an approximation of the solution $u(\x)$ of a boundary-value problem:
\begin{equation*}
    \begin{aligned}
       &R(\x,u ) =0\ , \ \forall \x \in \upOmega ;\\
       &B(\x ,   u) =0\ , \ \forall \x \in \mathrm{\mathrm{\Gamma}},
    \end{aligned}
\end{equation*}
where $R$ is the residual of a PDE, and $B$ is the boundary condition on $u(\x)$ we want to impose on $\rm \Gamma.$
With the vanilla PINN, we approximate the solution $u(\x)$ by a feedforward neural network $N_{\boldsymbol{\zeta}}(\x)$, and we tune the weights and biases $\boldsymbol{\zeta}$ by minimizing the loss function
\begin{align*}
L &= \alpha_1 \sum_{i=1}^{n_\upOmega} R (\x_\upOmega^i, N_{\boldsymbol{\zeta}}(\x_\upOmega^i)  )^2 + \alpha_2 \sum_{j=1}^{n_\mathrm{\Gamma}} B  (\x_\mathrm{\Gamma}^j, N_{\boldsymbol{\zeta}}(\x_\mathrm{\Gamma}^j)  )^2 + \alpha_3 \sum_{k=1}^{n_{\rm D}} ( N_{\boldsymbol{\zeta}}(\x_{\rm D}^k) - u(\x_{\rm D}^k)  )^2  \\[2mm] 
&= \alpha_1 \|  R (\x_\upOmega, N_{\boldsymbol{\zeta}}(\x_\upOmega)  ) \| + \alpha_2 \|B  (\x_\mathrm{\Gamma}, N_{\boldsymbol{\zeta}}(\x_\mathrm{\Gamma})  ) \| + \alpha_3 \|  N_{\boldsymbol{\zeta}}(\x_{\rm D}) - u(\x_{\rm D})   \|
\end{align*}
where $\{\x_\upOmega^i\}_{i=1}^{n_\upOmega}$ is a set of $n_\upOmega$ points that belong to $\upOmega$, $\{\x_\mathrm{\Gamma}^i\}_{i=1}^{n_\mathrm{\Gamma}}$ is a set of $n_\mathrm{\Gamma}$ points that belong to $\mathrm{\Gamma}$,  $\{\x_{\rm D}^i\}_{i=1}^{n_{\rm D}}$ is a set of $n_{\rm D}$ points where the solution is known (experimental measurements), $\alpha_i$ are weighting parameters that are user-defined, and $\|\cdot\|$ is the discrete $L2$ norm.

\subsection{Proposed PINN for the level-set advection problem}

In this work, we will limit our use of PINN to approximate the level-set function $F$, while the finite element method will give an approximation of the velocity field $\boldsymbol{u}$ and pressure $p$. 
\subsubsection{Time integration schemes for the advection equation}
Three approaches are considered in this work to handle time evolution. The first two are continuous-in-time, and the third approach is discrete-in-time. For the continuous-in-time approach, the time domain is decomposed into smaller time intervals with a time step size ${\rm \Delta} t$. At each time step $t_n$, we compute the FEM approximation of the velocity field $\boldsymbol{u}_n$. This velocity field is then used to obtain an approximation of the level-set function $F(\x,t)$, defined over the entire time interval $t \in [t_{n-1} ,t_n]$.

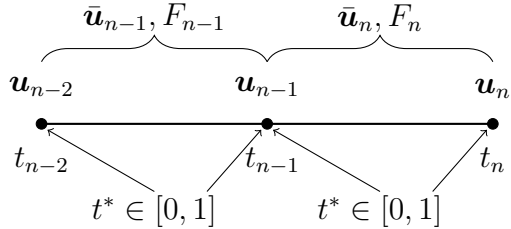
\begin{figure}
    \centering
\begin{tikzpicture}[decoration={brace,amplitude=10pt}]
  % Draw line segments
  \draw[thick] (0,0) -- (3,0) -- (6,0);
  
  % Points at t_{n-2}, t_{n-1}, t_n
  \filldraw (0,0) circle (2pt) node[below=5pt] {$t_{n-2}$};
  \filldraw (3,0) circle (2pt) node[below=5pt,xshift=2pt] {$t_{n-1}$};
  \filldraw (6,0) circle (2pt) node[below=5pt] {$t_n$};
  
  % Labels u_{n-2}, u_{n-1}, u_n
  \node[above] at (0,0.2) {$\boldsymbol{u}_{n-2}$};
  \node[above] at (3,0.2) {$\boldsymbol{u}_{n-1}$};
  \node[above] at (6,0.2) {$\boldsymbol{u}_n$};
  
  % Accolades (braces) for \bar{u}_{n-1} and \bar{u}_n
  \draw[decorate, yshift=2em] (0,0) -- node[above=0.6em] {$\bar{\boldsymbol{u}}_{n-1},F_{n-1}$} (3,0);
  \draw[decorate, yshift=2em] (3,0) -- node[above=0.6em] {$\bar{\boldsymbol{u}}_n, F_{n}$} (6,0);
  
  % Labels F_{n-1} and F_n below segments
  \node[below] at (1.5,-0.8) {$t^* \in [0,1]$};
  \draw[->] (1.5,-0.9) -- (0.08,-0.08);
  \draw[->] (2.2,-0.9) -- (3-0.08,-0.08);
  
  \node[below] at (4.5,-0.8) {$t^* \in [0,1]$};
  \draw[->] (4.5,-0.9) -- (3+0.08,-0.08);
  \draw[->] (5.2,-0.9) -- (6-0.08,-0.08);
\end{tikzpicture}

\caption{A representation of two time intervals, showing the discrete times $t_n$, the discrete velocities $\boldsymbol{u}_n$, the continuous in time interpolation of the velocities $\bar{\boldsymbol{u}}_n $, and the level-set function $F_n$ defined over each time interval.}
\label{fig:continous_rep}
\end{figure}

For a small time step size ${\rm \Delta} t$, the time interval becomes too small. Neural networks perform better during training if the inputs are normalized and have a magnitude of $O(1)$ \cite{imp_nor}, hence the need to normalize the time input.
 Therefore, we define a normalized time $t^* = \frac{t - t_{n-1}}{{\rm \Delta} t}$, which falls within the interval $t^* \in [0,1]$ when $t \in [t_{n-1},t_n]$. The advection equation of the level-set function becomes
\begin{equation}
 \partial_{t^*} F_n + {\rm \Delta} t \: \boldsymbol{\overline{u}}_n \cdot \nabla F_n =0, \ \text{where} \ 
 F_n(\x,t^*=0) = F_{n-1} (\x,t^*=1),
 \label{eq:advect}
\end{equation}
where $F_{n-1}$ is the level-set function defined over the previous time interval $[t_{n-2},t_{n-1}]$, $F_n$ is the discrete level-set function defined over the current time interval $[t_{n-1}, t_n]$, and $\boldsymbol{\overline{u} }_n$ is the continuous advection velocity defined for all $t \in [t_{n-1}, t_n]$, obtained by interpolating the discrete velocities calculated with the FEM $\boldsymbol{u}_{n}, \boldsymbol{u}_{n-1} , \dots \boldsymbol{u}_{0}$ in time. 
For example, we can define $\boldsymbol{\overline{u}}_n = \boldsymbol{u}_{n-1}(1-t^*) + \boldsymbol{u}_n t^*$ as the linear interpolation of the velocity fields.
Figure~\ref{fig:continous_rep} gives a graphical representation of two time intervals, showing the discrete times and velocities $t_n$ and $\boldsymbol{u}_n$, the level-set function $F_n$ defined on each interval, the interpolated velocity $\boldsymbol{\overline{u} }_n$, and the normalized time $t^*$ for each interval.

With the first continuous-in-time approach, which we will refer to as the weak PINN, the initial conditions are weakly imposed. The neural network will take the space-time coordinates $[x,y,z,t^*]$ as input, and the network is trained using the loss function 
\[ L = \| \partial_{t^*} F_n(\x,t^*) + {\rm \Delta} t \: \boldsymbol{\overline{u}}_n \cdot \nabla F_n(\x,t^*)  \| + \|F_n(\x,t^*=0) - F_{n-1}(\x,t^*=1) \|, \]
where $\boldsymbol{\overline{u}}_n = \boldsymbol{u}_{n-1}(1-t^*) + \boldsymbol{u}_n t^*$. We used a linear interpolation for $\overline{\boldsymbol{u}},$ but higher-order interpolations are possible.

With the second continuous-in-time approach, which we refer to as the strong PINN, we strongly impose the initial conditions as proposed in \cite{sukumar2022exact}, where the output of the neural network $N_{\boldsymbol{\zeta}}(\x,t^*)$ is multiplied by a distance function to the boundary $d(\x,t^*)$, and we add a lift function $V(\x,t^*)$ that satisfies the desired boundary condition. In our case, the distance function is $d(\x,t^*) = t^*$, 
$V(\x,t^*) = F_{n-1}(\x,t^*=1) \ (1-t^*)$, and the level-set function can be expressed as

\begin{equation}
    F_n(\x,t^*) = N_{\boldsymbol{\zeta}}(\x,t^*) \  t^* + F_{n-1}(\x,1) \ (1-t^*),
\label{eq:strong_time_imp}
\end{equation}
on $[t_{n-1},t_n].$
With this formulation, the level-set function at $t^*=0$ satisfies the initial condition $F_n(\x,0) = F_{n-1}(\x,1) ,$ and the loss function becomes
\[ L = \| \partial_{t^*} F_n(\x,t^*) + {\rm \Delta} t \ \boldsymbol{\overline{u}}_n \cdot \nabla F_n(\x,t^*)  \|. \]

The reason why we define $V(\x,t^*) = F_{n-1}(\x,1) \ (1-t^*)$ is to ensure that the level-set function $F_n$ only depends on the previous level-set function $F_{n-1}$. If we would consider $V(\x,t^*) = F_{n-1}(\x,1),$ for example, we would obtain that
\begin{align*}
    F_n(\x,t^*) &= N_{\boldsymbol{\zeta},n}(\x,t^*) \  t^* + F_{n-1}(\x,1)  \\ &= N_{\boldsymbol{\zeta},n}(\x,t^*) \  t^* + N_{\boldsymbol{\zeta},n-1}(\x,1) +  F_{n-2}(\x, 1)  \\ &=   N_{\boldsymbol{\zeta},n}(\x,t^*) \  t^* + N_{\boldsymbol{\zeta},n-1}(\x,1) + \ldots + N_{\boldsymbol{\zeta},1}(\x,1)+  F_{0}(\x, 1),
\end{align*}
where $N_{\boldsymbol{\zeta},n}$ is the neural network trained over the time interval $[t_{n-1},t_n].$
On the other hand, for $V(\x,t^*) = F_{n-1}(\x,1) \ (1-t^*)$, we obtain that
\begin{align*}
F_n(\x,t^*) &= N_{\boldsymbol{\zeta},n}(\x,t^*) \  t^* + F_{n-1}(\x,1)(1-t^*) \\
&= N_{\boldsymbol{\zeta},n}(\x,t^*) \  t^* + (N_{\boldsymbol{\zeta},n-1}(\x,1) \, (1) +  F_{n-2}(\x, 1) \, (1-1) ) (1-t^*) \\
&=N_{\boldsymbol{\zeta},n}(\x,t^*) \  t^* +N_{\boldsymbol{\zeta},n-1}(\x,1) (1-t^*).
\end{align*}
This results in less memory requirements, since we only need to keep the parameters of two neural networks in memory, $N_{\boldsymbol{\zeta},n}$ and $N_{\boldsymbol{\zeta},n-1}$, for each time interval, and the partial derivatives at $t^*=1$, which are used to approximate the capillary forces, will depend exclusively on the output of one neural network, ${\nabla F_n(\x,1)= \nabla N_{\boldsymbol{\zeta},n}(\x,1)}$. 

With the discrete-time approach, referred to as the finite difference PINN (FDPINN), we use a finite difference approximation of the time derivative. The neural network inputs are the spatial coordinates $[x,y,z]$. An implicit scheme is used, and we consider the advection velocity to be $\overline{\boldsymbol{u}} = \boldsymbol{u}_n.$ The loss function is then defined as
\[ L = \Big\| \frac{F_n(\x) - F_{n-1}(\x)}{{\rm \Delta} t} +  \boldsymbol{u_n} \cdot \nabla F_n(\x)  \Big\| . \]

% To address this issue, we propose mapping the input vector to a higher-dimensional space similar to the approach of \cite{fourier_net}.
\subsubsection{Neural network architecture}
Neural networks exhibit a spectral bias, where the associated optimization problem tends to converge faster toward the low-frequency components of the solution~\cite{spec_bias,spec_bias2}.
This spectral bias causes the problem to converge slowly when the fluid interface undergoes significant topological changes, such as bubble splitting, leading to prolonged computational times, and reduced accuracy. 
Multiple methods have been proposed to address the spectral bias of neural networks, most notably the Fourier Feature Neural Network \cite{fourier_net}, where the presence of high-frequencies in the output
is dealt with by mapping the input space to a high-frequency component space. However, this method did not produce satisfactory results for multi-fluid flow problems for two reasons. First, these methods work better for cases where the solution has a uniform frequency. This is not the case for multi-fluid flows, since the interface could be smooth in one part of the domain, and the problem can become stiff in another part of the domain. The second reason is that introducing high frequencies in the input vector will introduce high frequencies in the output vector, resulting in inaccurate free surface curvature approximations. 

In this work, we use a multi-level neural network, inspired by the work of Aldirany et al.~\cite{ziad}. We use two PINNs, giving $F_{n,1}$ and $F_{n,2}$ such that $$F_n = F_{n,1} + \beta F_{n,2},$$ where $\beta$ is a normalization parameter used to normalize the function $F_{n,2}.$
The first function $F_{n,1}(\x,t),$ approximates the level-set function across the entire domain. The input of the neural network is the space-time coordinates of points scattered over the entire domain. In this work, we choose the nodes of the FEM mesh as collocation points and they are assigned with random time coordinates $t^*$ within the interval $t^* \in [0,1]$.

The PINN is trained for the boundary-value problem 
\begin{equation}
\partial_{t^*} F_{n,1} + {\rm \Delta} t \  \boldsymbol{\overline{u}}_n \cdot \nabla F_{n,1} =0, \ \text{where} \ F_{n,1}(\x,0) = F_{n-1}(\x,1).
\label{eq:eq_F_1}
\end{equation} 
The second function $F_{n,2}$ computes a higher-quality approximation near the interface.  
Next, we define the input vector of the PINN that gives $F_{n,2},$ the boundary-value problem for $F_{n,2},$ and the normalization parameter $\beta.$

Since interfaces evolve in their normal direction, adding the interface normal vectors to the input of the PINN will help the network learn the evolution of the free surface. The normal vectors to the interfaces are approximated using $\nabla F_{n,1}$. Hence, the inputs of the second PINN, which produces $F_{n,2},$ are the space-time coordinates, and the gradient of the approximate solution $\nabla F_{n,1}$ at points around interfaces. These points are chosen such that $|F_{n,1}(\x,t^*)| < \epsilon,$ where $\epsilon$ is a positive constant that we choose to be $\epsilon=0.05.$
By replacing $F_n$ in equation \eqref{eq:advect}, we obtain
$$  
\partial_t^* (F_{n,1}+\beta F_{n,2}) + {\rm \Delta} t \bar{\boldsymbol{u}}_n \cdot \nabla (F_{n,1}+ \beta F_{n,2}) =0, \ \text{where} \ F_{n,1}(\x,0)+ \beta F_{n,2}(\x,0) = F_n(\x,1), $$

 and by isolating $F_{n,2}$, we obtain the boundary-value problem
\begin{equation}
\begin{aligned}
\partial_t^* \beta F_{n,2} +{\rm \Delta} t \bar{\boldsymbol{u}}_n \cdot \nabla \beta F_{n,2} &= - \partial_{t^*} F_{n,1} - {\rm \Delta} t \  \bar{\boldsymbol{u}}_n \cdot \nabla F_{n,1} 
\\& = -R(F_{n,1}),\ \text{where} \ \beta F_{n,2}(\x,0) = F_n(\x,1) -F_{n,1}(\x,0) ,
\label{eq:pdef2}
\end{aligned} 
\end{equation}
where $R(F_{n,1})$ is the residual of $F_{n,1}$.
\\
\\
N.B:\textit{ to accelerate the computation of the partial derivative of $F_{n,2}$ with respect to $\x$, we can apply the chain rule 
$$\frac{\partial F_{n,2}(\x)}{\partial x_i}= \frac{\partial F_{n,2} (\x, \nabla F_{n,1} ) }{\partial x_i }  + \sum_{j=1}^d \frac{\partial F_{n,2} (\x, \nabla F_{n,1} ) }{\partial ( \partial F_{n,1}(\x) / \partial x_j )} \frac{\partial ^2 F_{n,1}(\x)}{\partial x_i \partial x_j}, $$ and exclude $F_{n,1}$ from the computational graph.}
\\

As for the normalization of the function $F_{n,2}$ to obtain more accurate results, we need to approximate the value of $\beta$, so that $F_{n,2}$ is of the order of $O(1)$. We approximate $\beta$ using a finite difference scheme for equation \eqref{eq:pdef2}:
$$\frac{\beta (F_{n,2}(\x,1,\nabla F_{n,1}) - F_{n,2}(\x,0,\nabla F_{n,1}) )}{1-0}+ {\rm \Delta} t \: \bar{\boldsymbol{u}}_n \cdot \nabla \beta F_{n,2}(\x,0,\nabla F_{n,1}) = - \: R(F_{n,1}). $$
Considering that $F_{n,1}(\x,t)$ exactly satisfies the initial condition $F_n(\x,1) =F_{n,1}(\x,0)$, which is the case for the ``strong PINN'', we obtain from \eqref{eq:pdef2} that $F_{n,2}(\x,0,\nabla F_{n,1}) = 0$, leading to
$$\max(\beta F_{n,2}(\x,1,\nabla F_{n,1})) \approx \beta \approx \max (|R(F_{n,1})|).$$
The steps to obtain the level-set function approximation using the proposed method are:
\begin{itemize}
\item Compute the capillary forces using $F_{n-1}(\x,1)$;
\item Compute the velocity $\boldsymbol{u}_n$ using the Finite Element Method;
\item Interpolate the velocities to obtain $\bar{\boldsymbol{u}}_n$ over $t \in [t_{n-1} , t_n]$;
\item Determine the first approximation of the level-set function $F_{n,1}$ over the time interval $t^* \in [0,1]$ using the boundary value problem given by \eqref{eq:eq_F_1} using PINN;
\item Compute $\beta = \max( R(F_{n,1}))$;
\item Compute a better approximation around the interface by approximating the solution $F_{n,2}$ to the boundary value problem \eqref{eq:pdef2}, using a PINN with the input $[\x,\nabla F_{n,1}]$;
\item Compute the final approximation of the level-set function $F_n = F_{n,1} + \beta F_{n,2}.$
\end{itemize}

\subsection{Re-initialization methods for the level-set function}

The level-set function must be defined as a signed distance function to achieve accurate approximations of the normal vectors to the interfaces, of the curvature of the interfaces, and of the capillary force. This property is lost when the interface is advected, and due to numerical diffusion. Hence, the re-initialization of the level set function is necessary.

The typical approach for re-initialization involves finding an approximation of the Eikonal boundary-value problem:
\begin{alignat*}{2}
&\|\nabla \phi(\x) \| = 1,  \quad &\forall \x \in \upOmega; 
\\
&\phi(\x) = 0, \quad &\forall \x \in \mathrm{\Gamma}.
\end{alignat*}

The re-initialization methods fall mainly into two categories: 
\begin{itemize}
    \item The Pseudo-Time Evolution technique, where the level-set function is advected until steady-state is reached using:
\begin{align}
& \partial_\tau \phi + \text{sgn}(F)(\| \nabla \phi \| -1) =0; \\
& \phi(\x,\tau=0) = F(\x),
\end{align}
where $F(\x)$ is the level-set function that needs to be re-initialized, $\tau$ is a pseudo-time,  $\text{sgn}(F)=1$ if $F$ is positive, and $\text{sgn}(F) =-1$ if $F$ is negative.

However, re-initializing the level-set function discretely introduces errors in the position of the interfaces, leading to mass conservation errors. Several improvements have been introduced to alleviate this problem, such as the constrained reinitialization method \cite{Const_reini}, where the position of the interfaces is included in the formulation using the least-squares method.

\item The fast marching method is a technique that computes an approximate solution to the Eikonal equation by marching through the grid points, starting from the nearest points to the interface. It uses a priority queue to determine which grid points are processed next, based on the current solution estimates. The set of points is classified into three categories: the accepted, the tentative, and the distant points. The key idea is to approximate the solution for the tentative points using the approximate solution of the accepted points with an upwind scheme, and to move the points with the lowest value of $\phi$ from the tentative to the accepted category.
\end{itemize}

In this section, we consider two re-initialization methods. The first one aims at finding the distance from any point $\x_0$ in the domain to the interface by minimizing:
\begin{mini}|l|
	  {\x}{\| \x- \x_0 \|;}{}{}
	  \addConstraint{F(\x)=0.}
\end{mini}
Gradient-based optimization techniques are now feasible due to the continuous representation of $F(\x)$, and by the fact that we have access to $\nabla F(\x)$ by back-propagation and automatic differentiation. The optimization is performed based on two formulations of the problem. The first formulation, which we will refer to as the penalty reinitialization, uses the penalty method
\begin{mini}|l|
	  {\x}{\| \x-\x_0\| + \lambda \|F(\x)\|,}{}{}
\end{mini}
where $\lambda$ is a penalization parameter. 
The second formulation, which we will refer to as the unconstrained reinitialization, aims at removing the conflict between the terms $ \| \x-\x_0\| $ and  $\|F(\x)\|$, where we rather seek to align the vector $\x-\x_0 $ with the normal vector $\nabla F$, and hence we minimize 
% $$\|F(\x)\| + \big|1 - \frac{\nabla F}{\|\nabla F\|} \cdot \frac{\x-\x_0}{\|\x-\x_0 \|} \|.$$ 
\begin{mini}|l|
	  {\x}{\|F(\x)\| + \big\|1 - \frac{\nabla F}{\|\nabla F\|} \cdot \frac{\x-\x_0}{\|\x-\x_0 \|} \big\|}{}{}.
\end{mini}

To accelerate the convergence, we can take advantage of the fact that the position of the interfaces does not change significantly from one time step to the next, which makes the points $\x$ obtained from the previous reinitialization, a good initial guess for the current reinitialization.

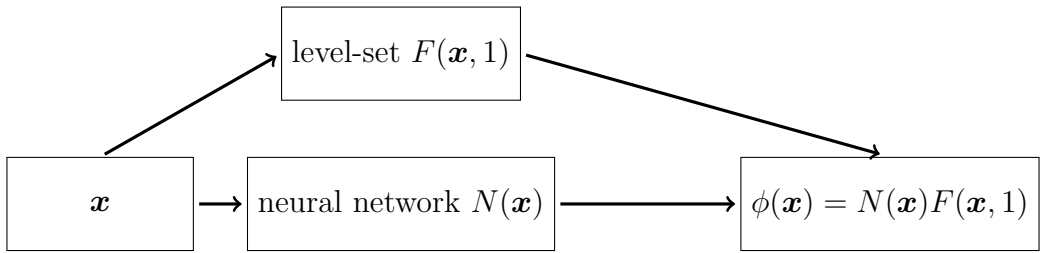
\begin{figure}[ht!]
\centering
\tikzstyle{block} = [draw, fill=white, rectangle, 
    minimum height=3em, minimum width=6em]
\tikzstyle{sum} = [draw, fill=white, circle, node distance=1cm]
\tikzstyle{input} = [coordinate]
\tikzstyle{output} = [coordinate]
\tikzstyle{pinstyle} = [pin edge={to-,thin,black}]

\begin{tikzpicture}[auto, node distance=2cm]

    \node [block, name=input] (controller) {$\x$};
    
    \node [block, right of=controller,
            node distance=4cm] (system) {neural network $N(\x)$};
            
    \node [block, above of=system] (boundary) {level-set $F(\x,1)$};
    %\node [block, below of=system] (distance) {$t$};
    
    \node [block, right of=system,node distance=6.5cm] (sol) {$\phi(\x)=N(\x)F(\x,1)$};

    \draw [->,very thick,shorten >=2pt, shorten <=2pt] (controller.north) -- node[name=u] {} (boundary.west);
    \draw [->,very thick,shorten >=2pt, shorten <=2pt] (controller) -- node[name=u] {} (system);
    
    \draw [->,very thick,shorten >=2pt, shorten <=2pt] (system) -- node[name=u] {} (sol);
    \draw [->,very thick,shorten >=4pt, shorten <=2pt] (boundary.east) -- node[name=u] {} (sol.north);
  
\end{tikzpicture}
\caption{The architecture of the PINN-R used to reinitialize the level-set function.}
\label{fig:eikonalpinn}
\end{figure}
The third approach, which we will refer to as the PINN Reinitialization (PINN-R), aims at approximating the solution to the Eikonal equation,$$\| \nabla \phi \| = 1,$$ using a PINN. The position of the interface is conserved by multiplying the output of the PINN, which is designed to give strictly positive outputs, by the level-set function $F(\x,1)$. This is shown in figure~\ref{fig:eikonalpinn}, where $F(\x,1)$ is the level-set function that needs to be reinitialized, the neural network $N(\x)$ gives a strictly positive scalar output, and the reinitialized level-set function $\phi(\x)$.

\section{Numerical results}

\subsection{Comparison of the time-integration schemes}
We first investigate which of the previously described time-integration schemes, the weak PINN, the FDPINN, and the strong PINN, is more accurate, using a manufactured problem.
The numerical test involves a cylinder with a radius of $R=1$, with height $h=1$. Inside the cylinder, a bubble with a radius of $r=0.25$, is positioned at the center of the cylinder, $\x_{\rm c}=[0,0,0.5]$. The bubble is subjected to the velocity field $\boldsymbol{u}=[-0.5x,-0.5y,z]$.
With this velocity field, the exact position of any particle, obtained by solving an ODE for each direction:
\begin{align*}
    &\frac{dx}{dt}=-0.5x; \\
    &\frac{dy}{dt}=-0.5y; \\
    &\frac{dz}{dt}=z,
\end{align*}
is given by
$$\boldsymbol{x}(t) = [x_0 e^{-0.5t} ,y_0 e^{-0.5t} ,z_0 e^{t}],$$ where $x_0$, $y_0$ and $z_0$ represent the particle's initial position.
The accuracy of the solution is measured using $E = \frac{1}{{ n}}\sum_{i=1}^{ n} |F(\x_i) |,$ where $\x_{i}$ is a set of ${n}$ points that belong to the interface at time $t=0$.
The chosen time step size is ${\rm \Delta} t=0.1,$ for a total simulation time of $t=1$. The number of collocation points is 70,000, and we train the PINN for 4000 iterations using the Adam algorithm, with a learning rate of $\alpha = 10^{-4}.$ We do not use the multi-level approach, and use instead a feedforward neural network with 5 layers and 50 neurons per layer.

\begin{figure}[ht!]
\begin{subfigure}[t]{0.3\textwidth}
\centering
    \includegraphics[width=\textwidth, height=\textwidth]{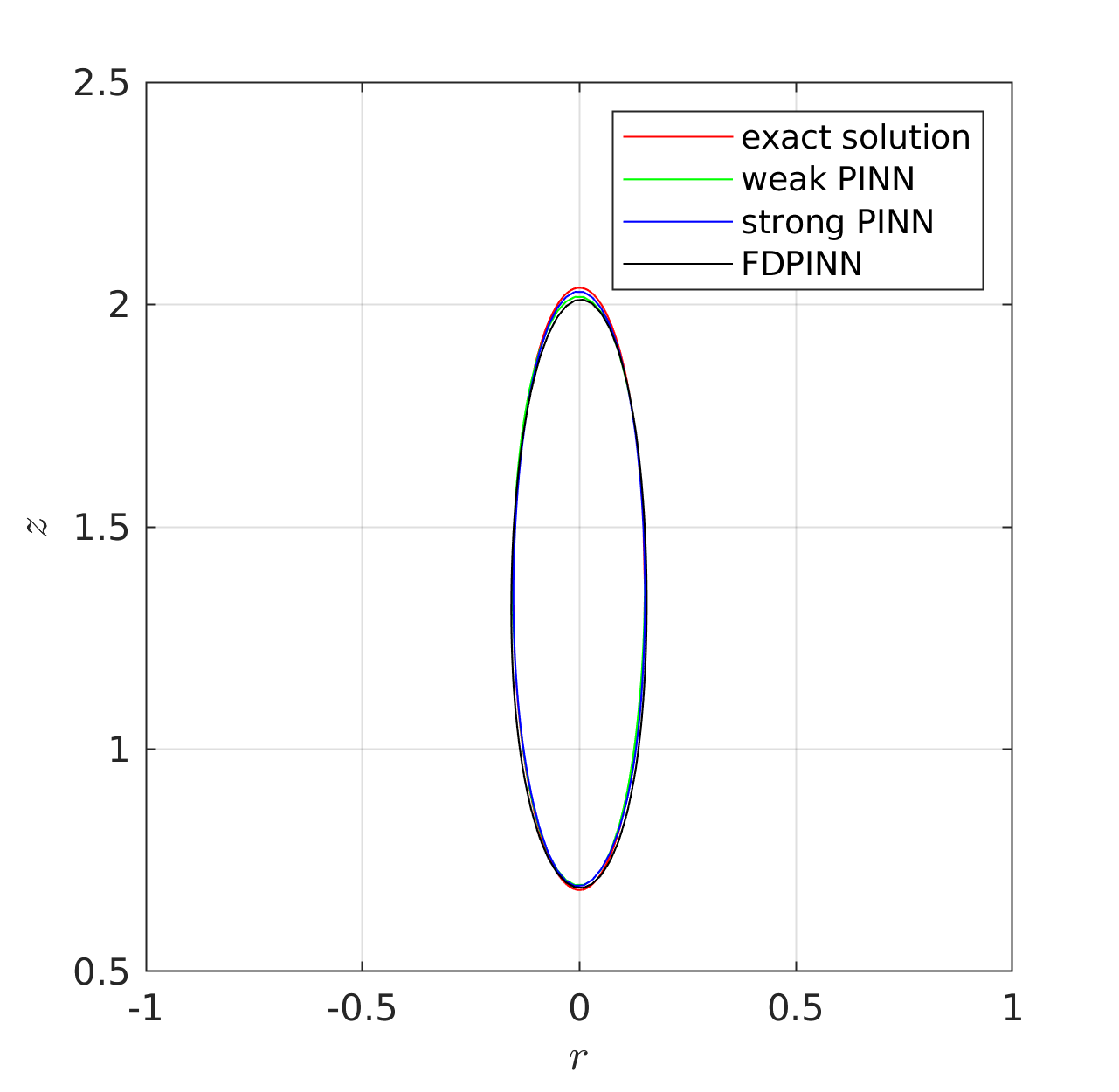}
    \label{fig:manu_sol_a}
    \caption{}
\end{subfigure}
\begin{subfigure}[t]{0.3\textwidth}
\centering
    \includegraphics[width=\textwidth, height=\textwidth]{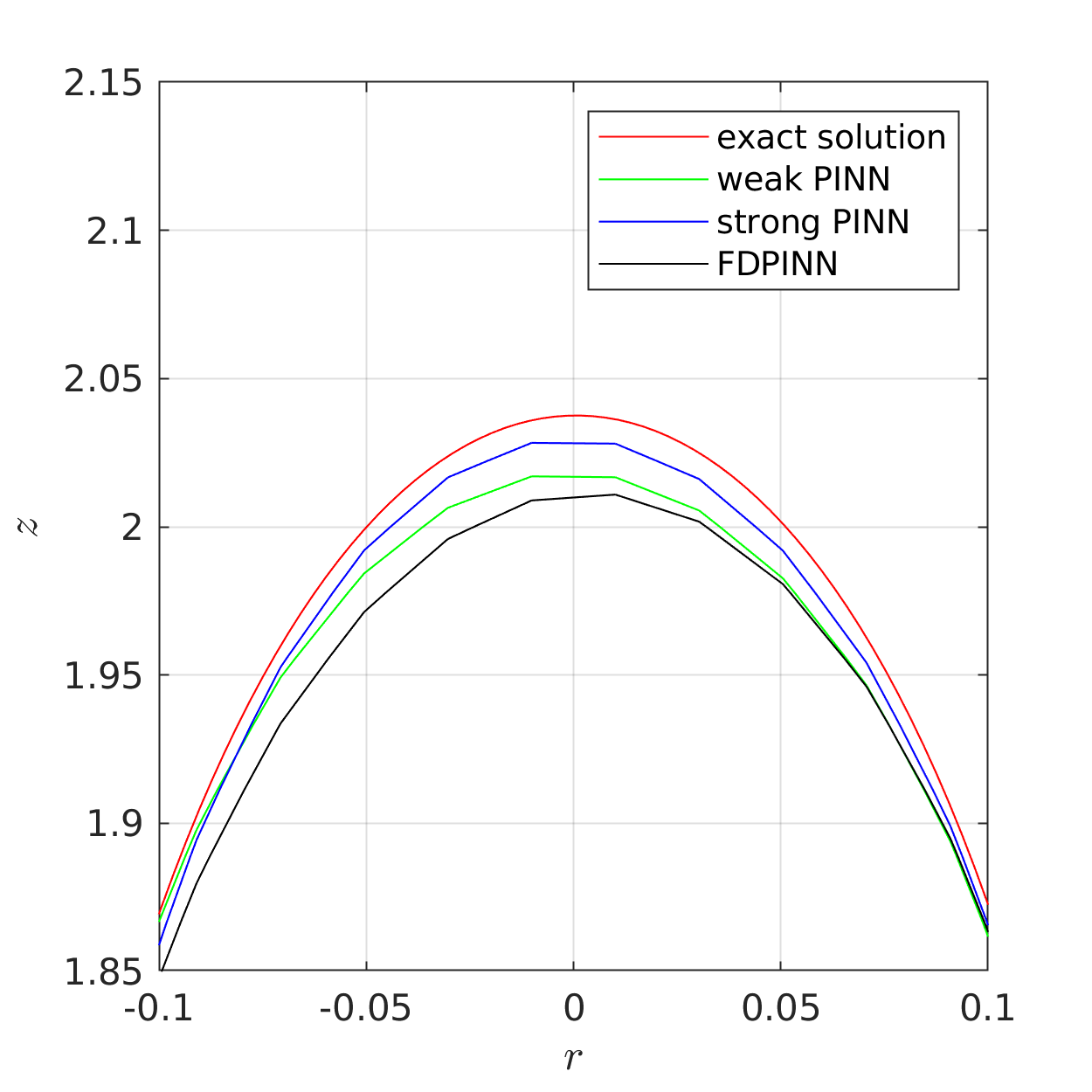}
    \label{fig:manu_sol_b}
    \caption{}
\end{subfigure}
\begin{subfigure}[t]{0.3\textwidth}
\centering
    \includegraphics[width=\textwidth, height=\textwidth]{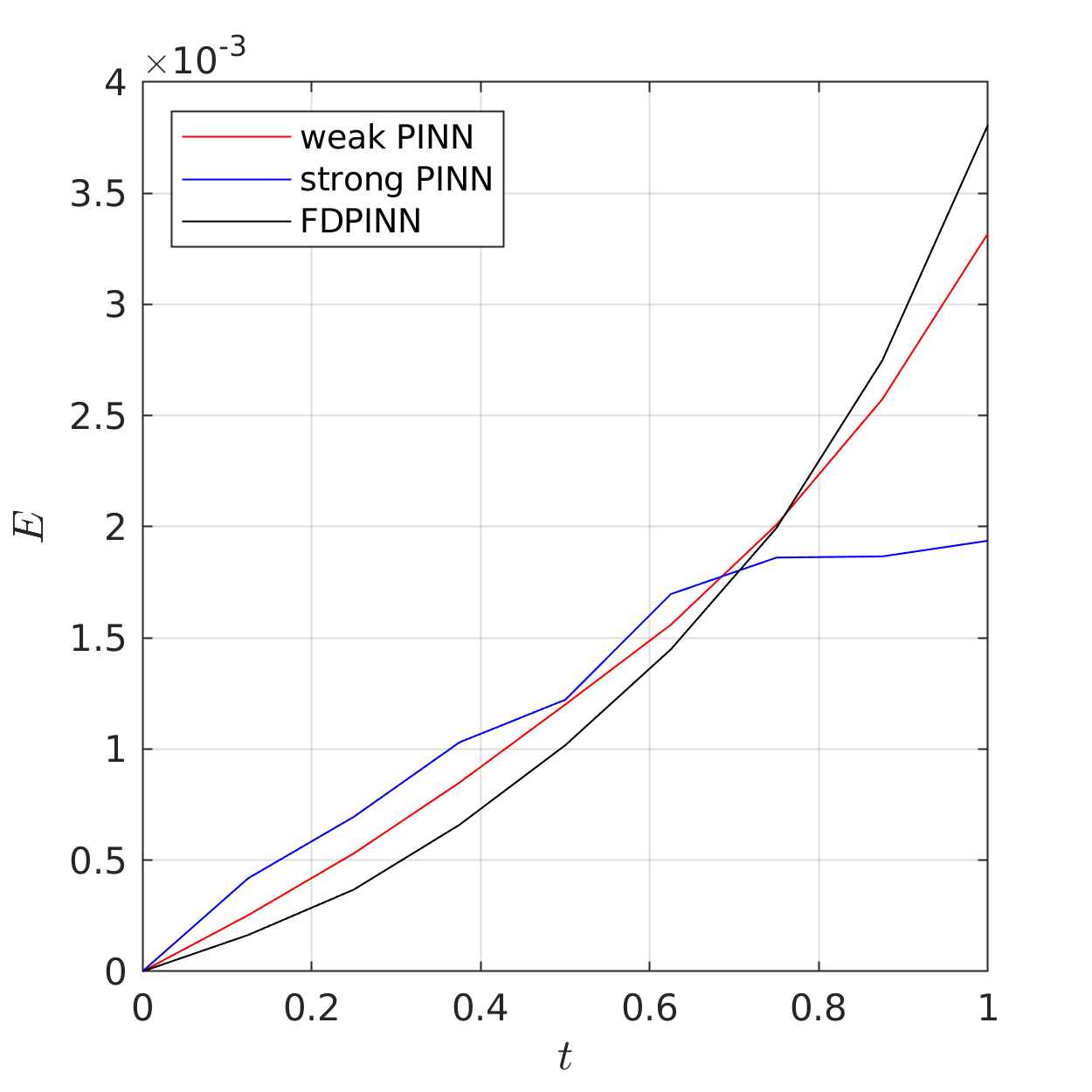}
    \label{fig:manu_sol_c}
    \caption{}
\end{subfigure}
    \caption{(a) The interface shape obtained at $t=1$ for the comparison of the three time-integration schemes, (b) a zoom-in on the interface shape at $t=1$,  (c) the error on the position of the interface, for the problem with a manufactured solution using the three different time integration schemes. }
    \label{fig:manu_sol}
\end{figure}

 As shown in figure \ref{fig:manu_sol}, strongly imposing the initial conditions yields the most accurate results in comparison with the other two formulations. This formulation does not have a time discretization error, contrary to the FDPINN, and it does not suffer from errors on the initial conditions, contrary to the weak PINN.

% \subsection{The effect of hyper-parameters}
% Since the second formulation yielded the most promising results, we will test the impact of the number of collocation points, the choice of ${\rm \Delta} t$, and the number of training iterations on the accuracy of the solution for the strong PINN formulation. As we can see in Figure~\ref{fig:er_param}, only the number of training iterations affects the quality of our result. The value of ${\rm \Delta} t$ plays no significant role, and the number of collocation points has no impact when enough points are present.

% \begin{figure}[ht!]
%     \includegraphics[width=7.5cm]{er_coloc_dt.png}
%     \includegraphics[width=7.5cm]{er_coloc.png}
%     \includegraphics[width=7.5cm]{er_iter.png}
%     \caption{ The error on the position of the interface for the second formulation with different time steps, collocation points, and training iterations. The first graph on the top left shows The error on the interface position for different ${\rm \Delta} t$ with 700000 collocation points and 4000 training iterations. the second graph on the top right shows The error on the interface position for different numbers of collocation points, ${\rm \Delta} t=0.1$, and 4000 training iterations. The third graph shows the error on the interface position for different numbers of training iterations with ${\rm \Delta} t=0.1$ and 70000 collocation point  }
%     \label{fig:er_param}
% \end{figure}

\subsection{Comparison between the multi-level PINN and the strong PINN  }
In this section, we compare the proposed multi-level PINN with a strong imposition of the initial condition,  with the single-level strong PINN. We use the verification problem of subsection 3.4.1. The results are obtained with $70,000$ collocation points, and $ {\rm \Delta} t = 0.1.$ We train the multi-level PINN with ${n}=6000$ total iterations, ${ n}_1 = 3000$ for the first function $F_{n,1}$ and ${n}_2=3000$ for the second function $F_{n,2}.$ We train the single-level PINN with ${n}=6000$ iterations. 

\begin{figure}[ht!]
\begin{subfigure}[t]{0.3\textwidth}
\centering
    \includegraphics[width=\textwidth,height=\textwidth]{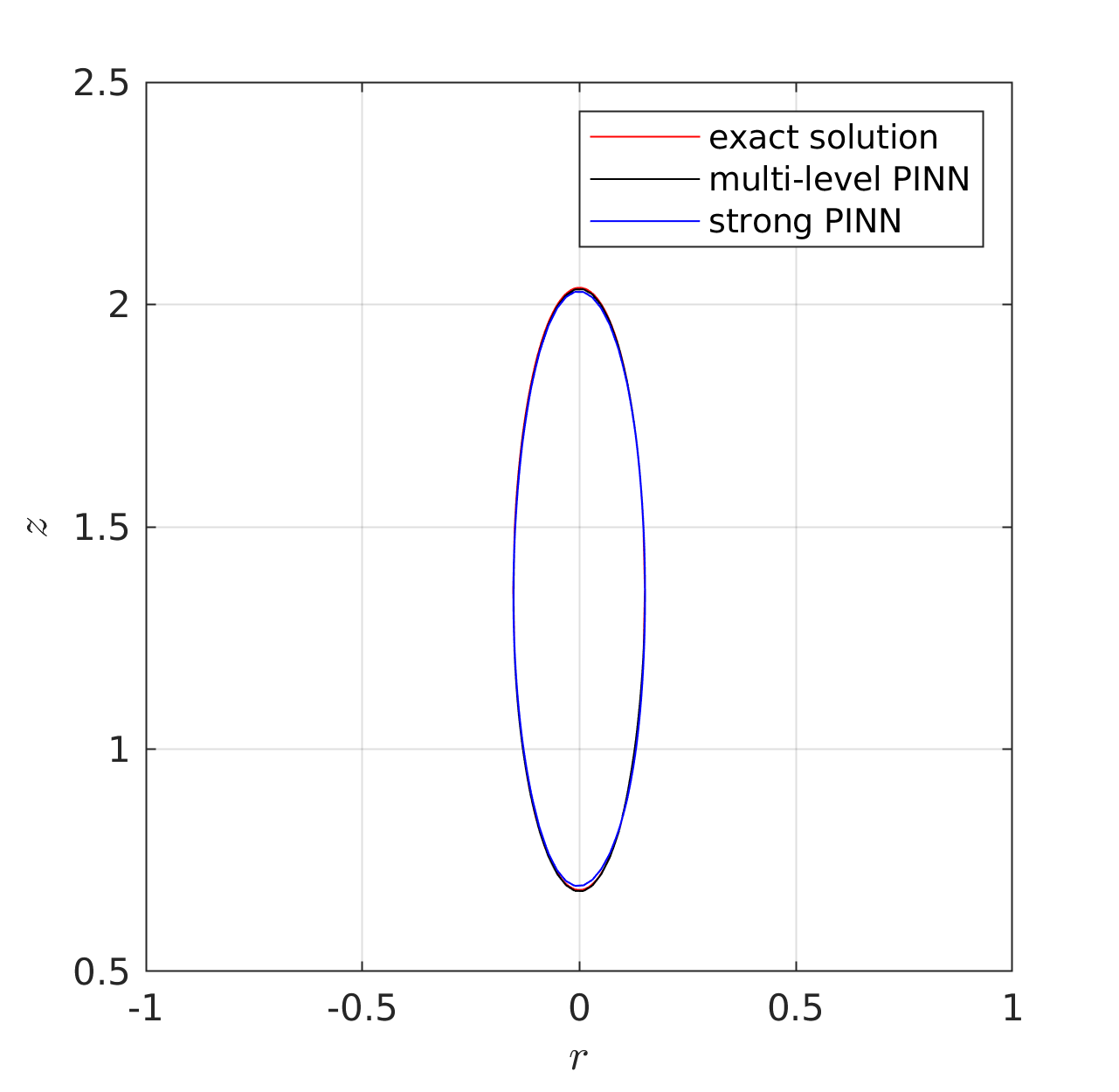}
    \caption{}
\end{subfigure}
\begin{subfigure}[t]{0.3\textwidth}
\centering
    \includegraphics[width=\textwidth,height=\textwidth]{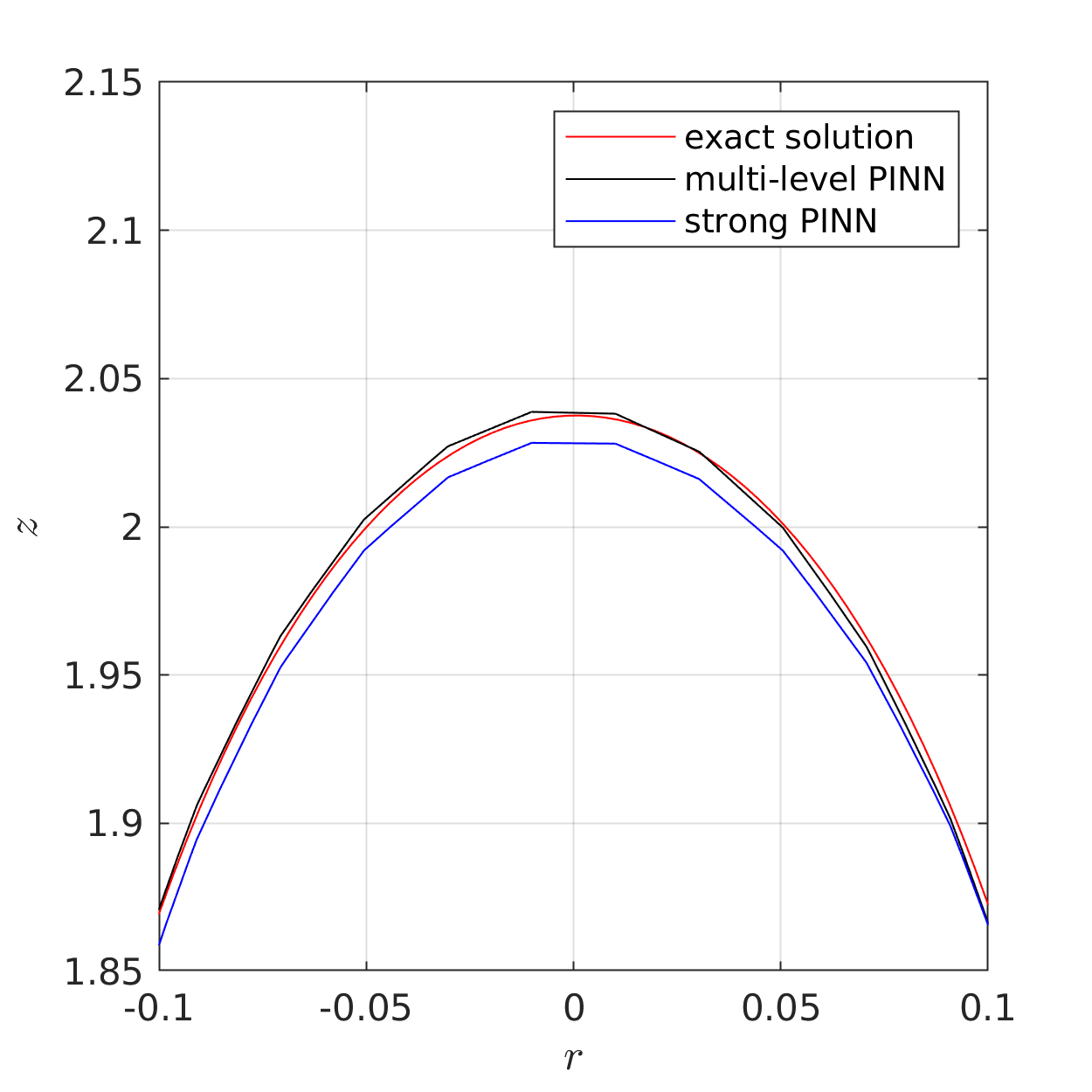}
    \caption{}
\end{subfigure}
\begin{subfigure}[t]{0.3\textwidth}
\centering
    \includegraphics[width=\textwidth]{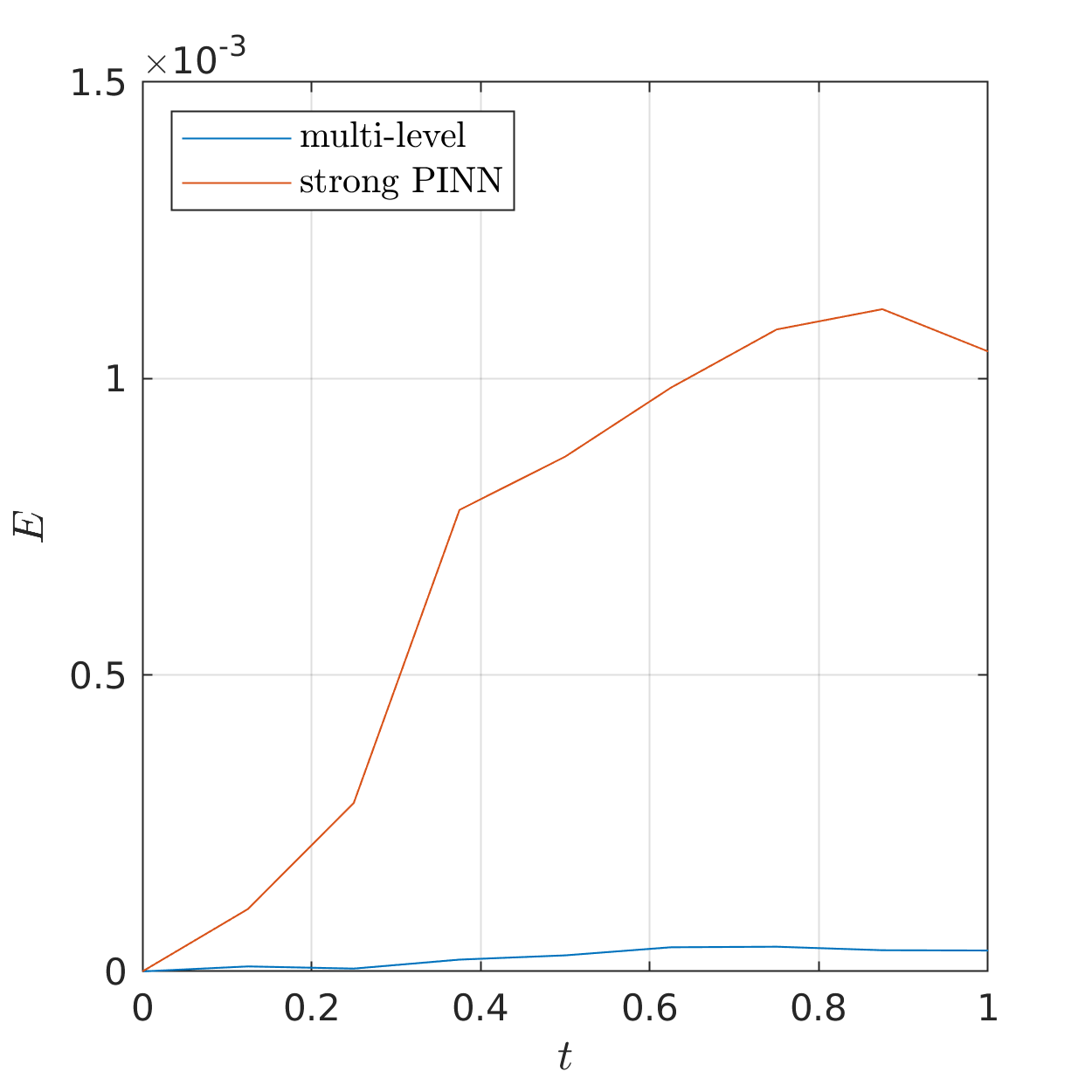}
    \caption{}
\end{subfigure}
    \caption{ (a) The interface shape obtained at time $t=1$ with the strong single-level PINN, and the proposed multi-level PINN, (b) a zoom-in on the interface shapes, (c) the error $E$ on the interface for the proposed methods, for the problem with a manufactured solution. }
    \label{fig:multi_vs_1}
\end{figure}

Figure~\ref{fig:multi_vs_1} illustrates the capabilities of the proposed multi-level PINN, for the same computational cost of the standard PINN formulation, where the error, measured as ${E = \frac{1}{ n}\sum_{i=1}^{ n} |F(\x_i) |}$, is equal to ${E =5\times 10^{-5}}$ for the multi-level PINN, and $E =10^{-3}$ for the regular PINN.
 
\subsection{Re-initialization methods}
\begin{figure}
\centering
    \includegraphics[width=7.5cm]{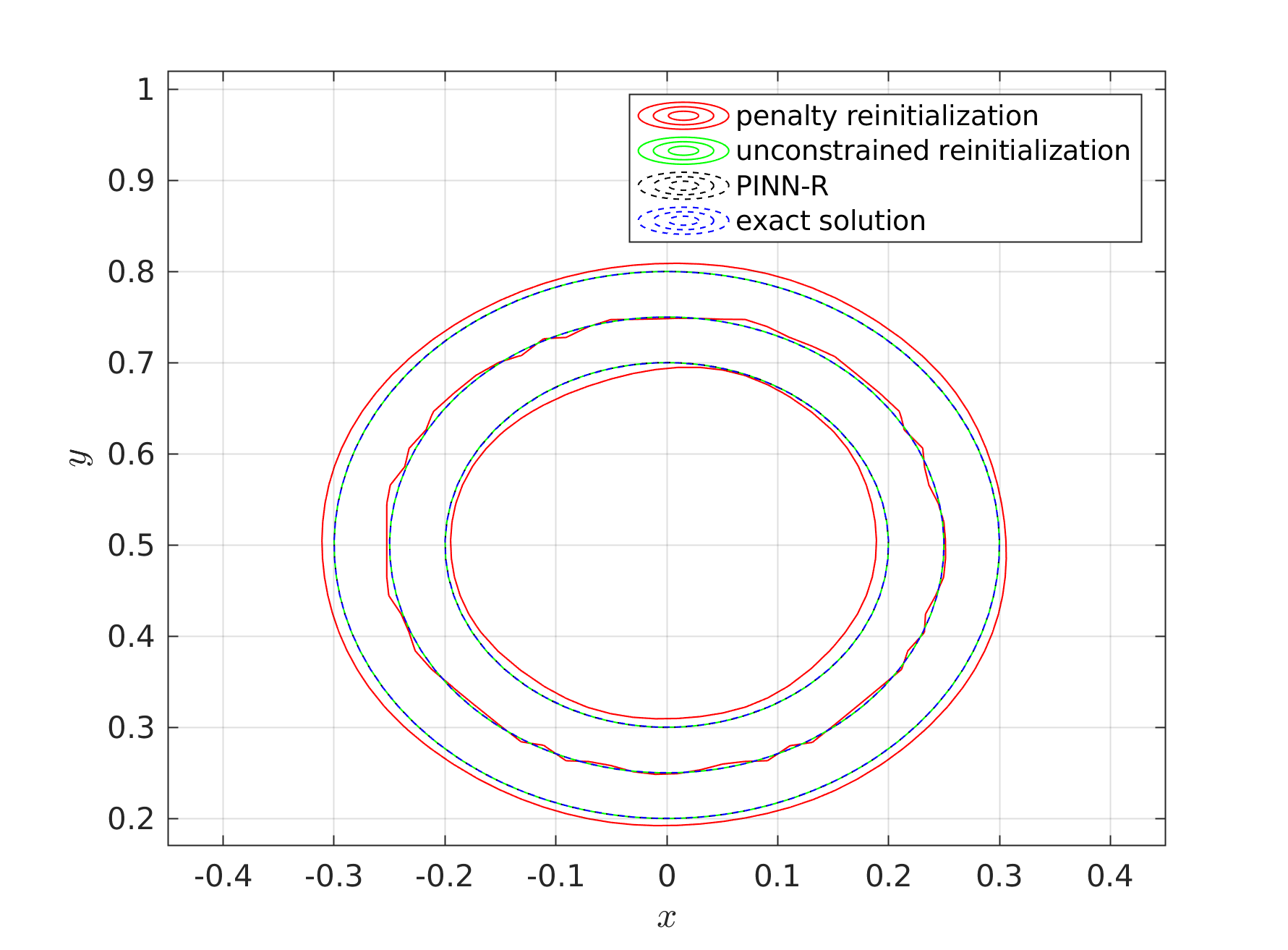}

    \caption{ The level-set function after the reinitialing $F = (x+1)^2(r- 0.25)^{3/2} $, using the penalty reinitialization, the unconstrained reinitialization, and the PINN-R method. }
    \label{fig:reini}
\end{figure}

\begin{figure}[ht!]
\begin{subfigure}[t]{0.5\textwidth}
\centering
    \includegraphics[width=7.5cm]{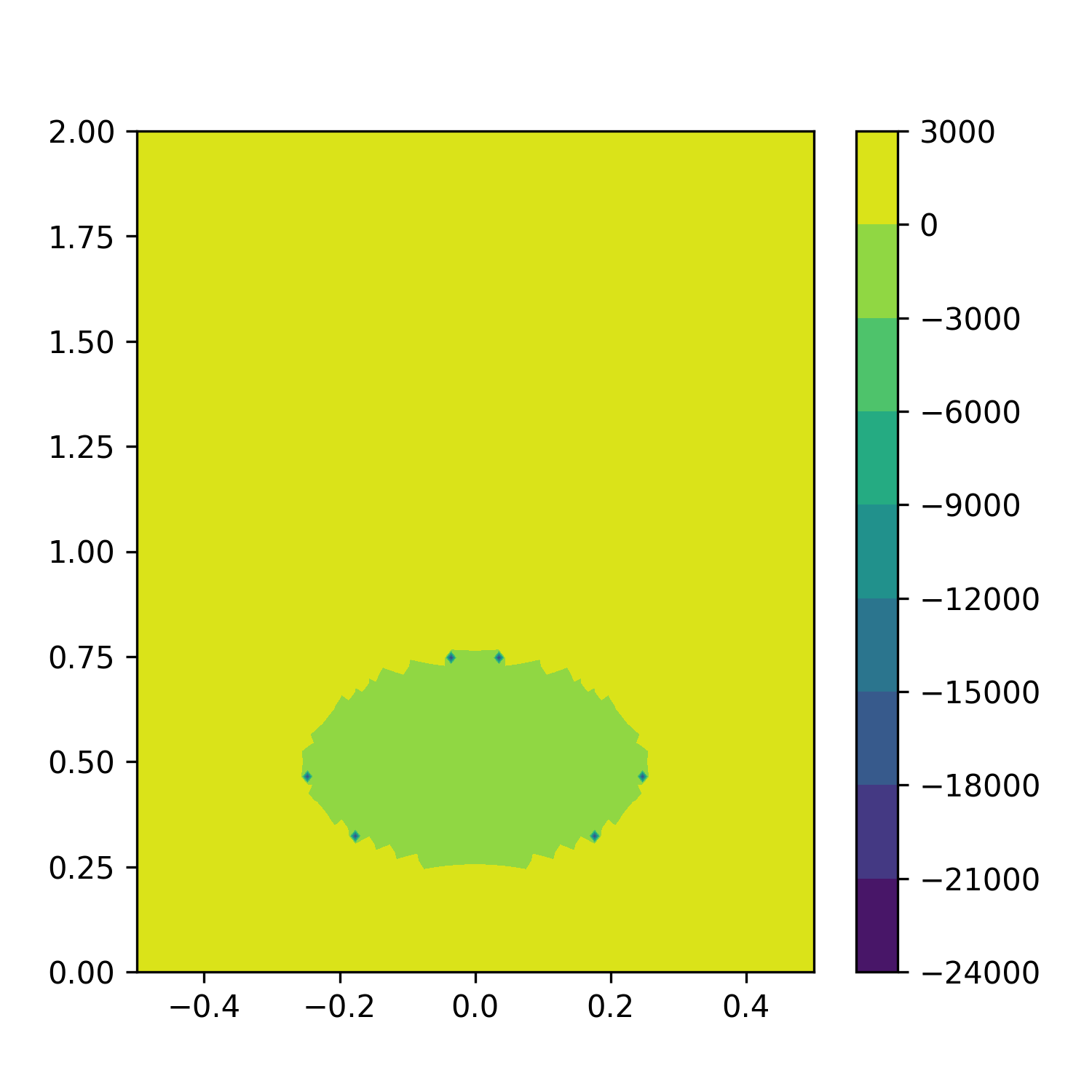}
\caption{}
\end{subfigure}
\begin{subfigure}[t]{0.5\textwidth}
\centering
    \includegraphics[width=7.5cm]{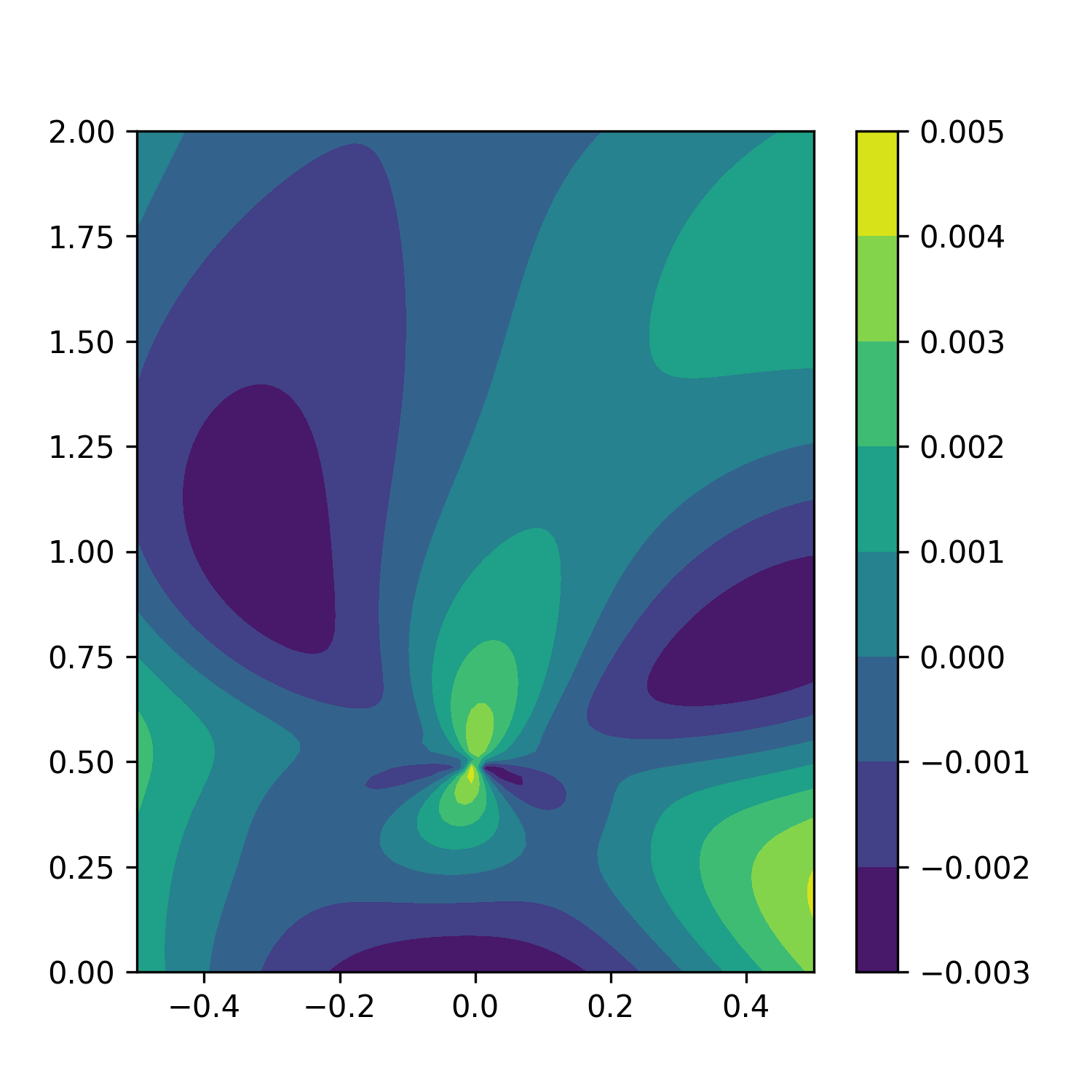}
    \caption{}
\end{subfigure}
    \caption{ The error on the curvature of the level-set function after reinitializing the level-set function $F = (x+1)^2(r- 0.25)^{3/2}, $ with (a) the unconstrained reinitialization method, (b) the PINN-R method.  }
    \label{fig:reini_cu}
\end{figure}

We apply the three reinitialization methods to a problem with a manufactured solution to assess their effectiveness. The domain is a $2$D unit square. The initial level-set function is given by $$F = (x+1)^2(r- 0.25)^{3/2}, $$ where $r$ is the distance to the point $\x_c = [0.5,0.5].$  The exact solution, after the reinitialization, is $\phi = r-0.25 .$
We compare the contours of the reinitialized level-set function obtained by the three approaches, as shown in figure \ref{fig:reini}. We notice that the results obtained by the penalty reinitialization are not accurate, when compared to the other reinitialization methods. Next, we compare the curvature obtained by the unconstrained reinitialization and the PINN-R. As seen in figure~\ref{fig:reini_cu}, reinitializing the level-set function with the PINN-R results in a smoother approximation, when compared to optimization-based methods. 

The numerical experiment demonstrates the superiority of the PINN-R for several key reasons.
 The optimization-based methods have more degrees of freedom, equivalent to the number of collocation points. In contrast, the PINN-R method involves a significantly lower number of degrees of freedom, which are the weights and biases of the neural network. Consequently, the PINN-R method is more computationally efficient.
The PINN-R method conserves the interface location during reinitialization, which leads to mass conservation, which is not inherently guaranteed by the other methods. 
The PINN-R approach leads to better curvature and capillary force approximations.
\\

\subsection{The steady-state Laplace problem}
For this problem, a circular bubble with a radius of $r=0.2$ is placed in a zero-velocity field without any external force. 
 The viscosity, density, and surface tension coefficient of the fluids are $\mu_1 = 10$, $\mu_2 = 0.1$,  $\rho_1 = 1000$, $\rho_2 = 1$, and $\gamma=1.96$.
The bubble should remain circular, and the magnitude of the velocity field should remain zero, with a pressure jump across the interface.
However, a velocity field appears due to numerical errors, usually referred to as spurious currents~\cite{zahedi2012spurious,abadie2015combined}. 
The capillary force is balanced out by the pressure jump across the interface, as pointed out by the Navier-Stokes equations when $\boldsymbol{u} = \boldsymbol{0}$, $$\nabla p = \kappa \gamma \delta \boldsymbol{n}.  $$
Hence, the pressure jump across the interface is given by ${\rm \Delta} p =  \kappa \gamma$. The spurious currents appear when the numerical pressure gradient cannot counteract the capillary force. Other factors are the errors in the curvature of the interface approximation, errors in the normal vectors to the interface, or a mismatch between the vector space of the pressure gradient and capillary force in the FEM formulation~\cite{zahedi2012spurious}. 

In the first test, the elements size is $h = 0.02,$ and the velocity field and pressure will be discretized using the $\rm P2/ \rm P1$ element combination. The exact level-set function, given by 
\begin{equation}
    F=\sqrt{(x-0.5)^2 + (y-0.5)^2}-0.2, 
\label{eq:ex_phi}
\end{equation}
is used. The capillary force is modeled using the continuum surface force model, and the partial derivatives of the level-set function are computed using automatic differentiation. The capillary force is interpolated using $\rm P0$ and $\rm P1$ elements to test the influence of the choice of the vector space on the spurious currents.

\begin{figure}[ht!]
\begin{subfigure}[t]{0.5\textwidth}
\centering
    \includegraphics[width=\textwidth]{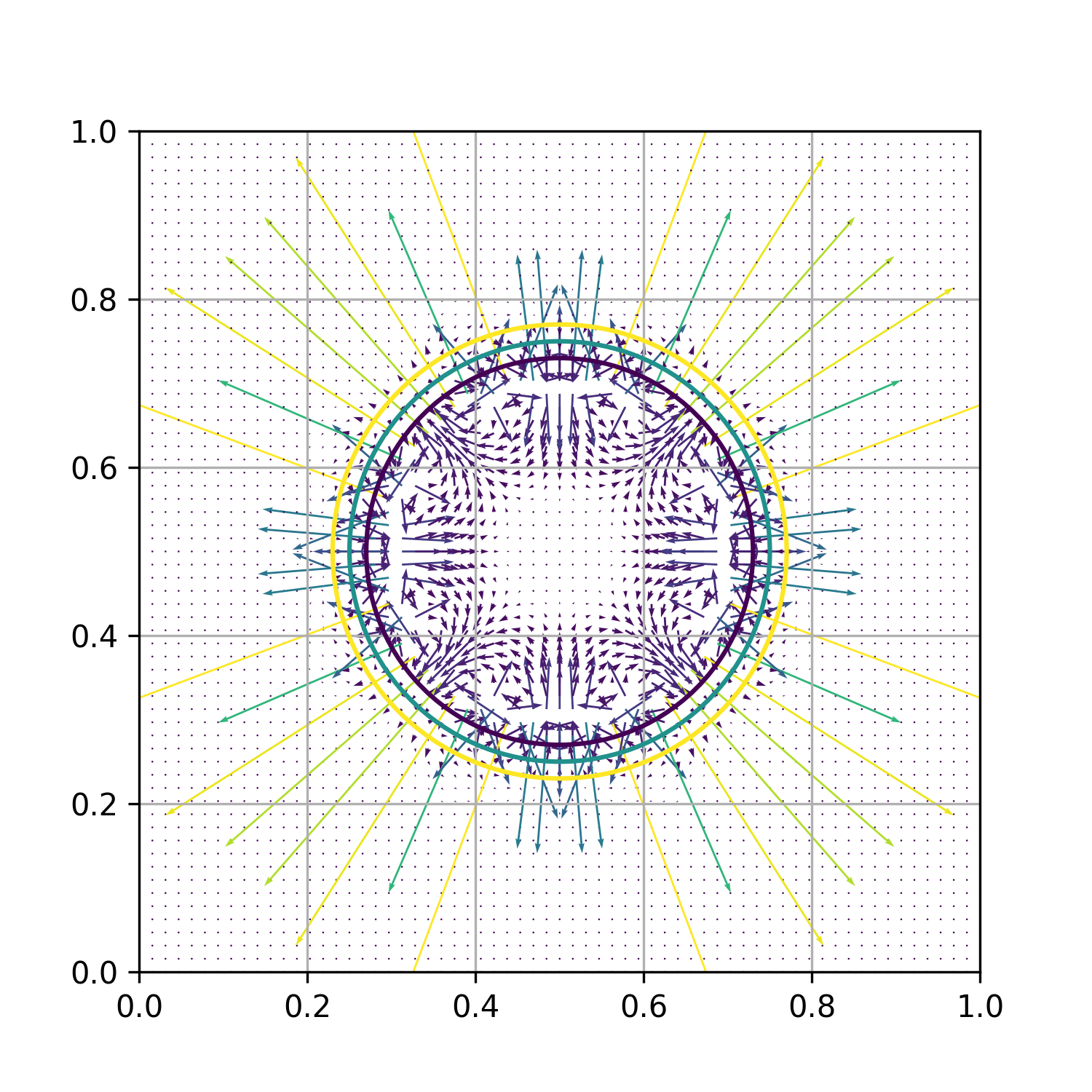}
\caption{}
\end{subfigure}
\begin{subfigure}[t]{0.5\textwidth}
\centering
    \includegraphics[width=\textwidth]{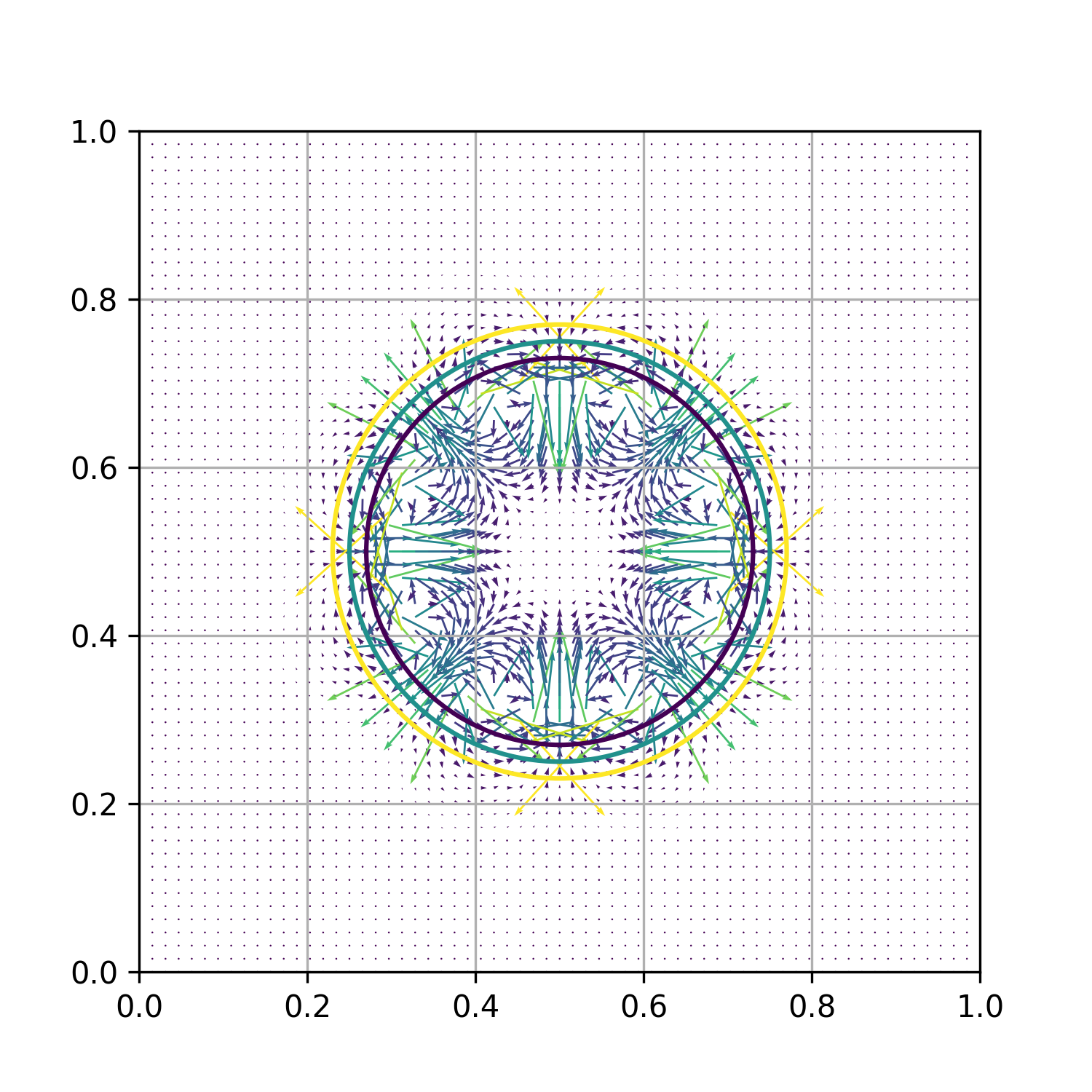}
\caption{}
\end{subfigure}
    \caption{ The spurious velocity field for the steady state Laplace problem obtained by interpolating the capillary force (a) using the $\rm P_1$ element, and (b) using the $\rm P_0$ element.}
\end{figure}
The L2 norm of the spurious currents is $\| \boldsymbol{u}\|=10^{-4}$ for the $\rm P_0$ approximation and $\| \boldsymbol{u}\|=5\times 10^{-4}$ for the $\rm P_1$ approximation.
In both cases, the pressure jump across the interface corresponds to the theoretical value with a relative error of $0.3 \% $. 

% \begin{figure}
% \centering
%     \includegraphics[width=5.5cm]{p.png}
% \caption{ The pressure with the $\rm P_1$ and the $\rm P_0$ discretization of the capillary force, compared to the exact solution.}
% \end{figure}
These results show that the choice of the elements used for the capillary force interpolation is important for reducing the spurious currents magnitude, even when the level-set function is exact.

\begin{figure}[ht!]
\begin{subfigure}[t]{0.5\textwidth}
\centering
    \includegraphics[width=\textwidth]{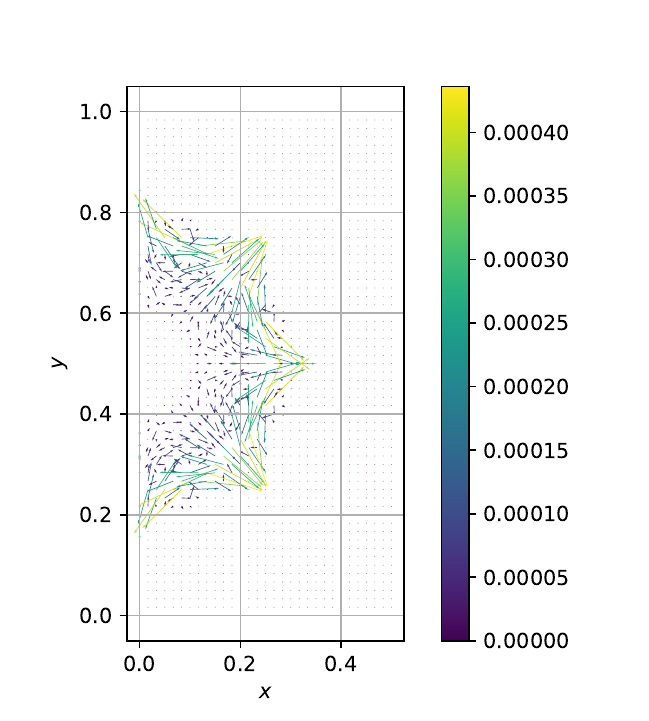}
\caption{}
\end{subfigure}
\begin{subfigure}[t]{0.5\textwidth}
    \includegraphics[width=\textwidth]{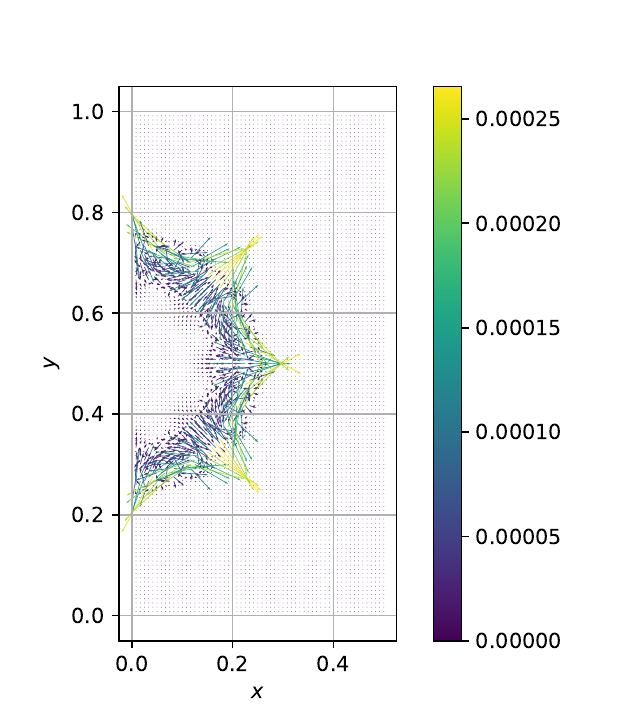}
\caption{}
\end{subfigure}
\begin{subfigure}[t]{0.5\textwidth}
    \includegraphics[width=\textwidth]{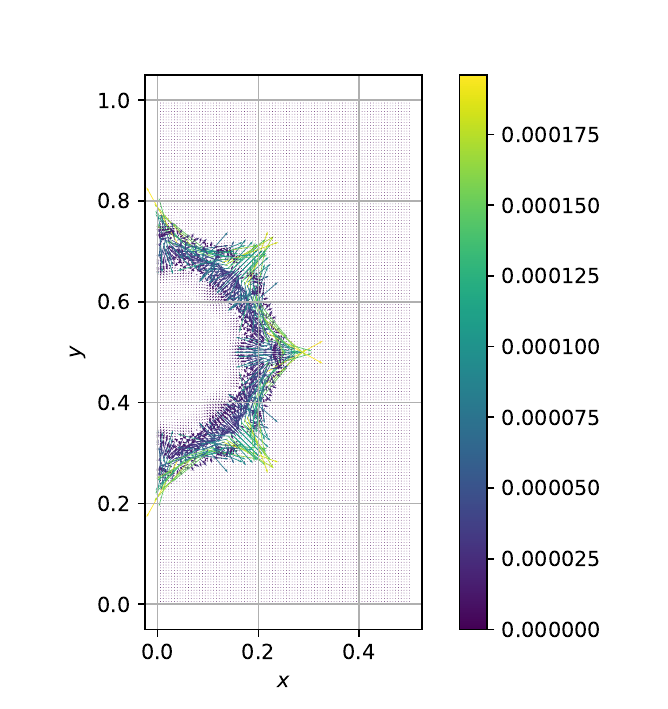}
\caption{}
\end{subfigure}
\begin{subfigure}[t]{0.5\textwidth}
    \includegraphics[width=\textwidth,height=\textwidth,]{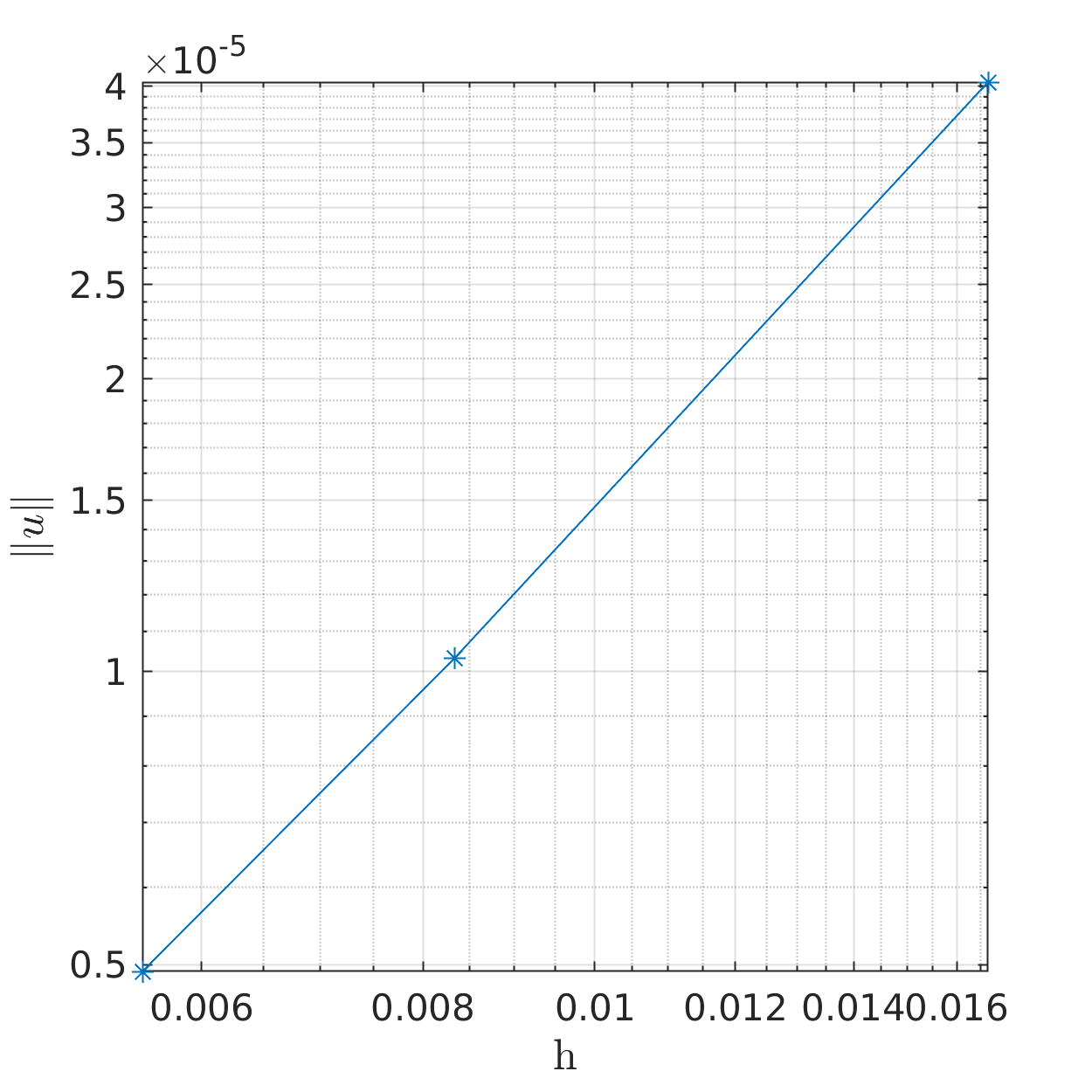}
\caption{}
\end{subfigure}
\caption{ The spurious velocities in the steady-state Laplace problem for different elements size, (a) with elements size $h = 1/60 ,$ (b) $h = 1/120 ,$ (c) $h = 1/180 ,$ and in (d) we show the $L2$ norm of the velocities. }
\label{fig:Vmesh}
\end{figure}

Next, we study the effect of mesh refinement on the magnitude of the spurious velocity. The velocity field and pressure are approximated using The $\rm P_2/ \rm P_1$ element combination. The exact level-set function \eqref{eq:ex_phi} is used.
 The capillary force is modeled using the continuum surface force model, where the partial derivatives of the level-set function are computed using automatic differentiation, and the surface tension is interpolated using the $\rm P0$ element. As seen in figure \ref{fig:Vmesh}, the element size plays an important role in reducing the spurious currents.

Finally, we test the effect of the PINN approximation error on the spurious currents. We train the neural network to approximate the initial condition with the loss function $$ L = \| F_0 - (\sqrt{(x-0.5)^2 + (y-0.5)^2}-0.2) \|, $$ using a different number of training iterations. The element size is $ h = 1/180$, and the capillary force is interpolated using the $\rm P_0$ element. The number of training iterations are $n_1 = 1000, \ 2000 ,$ and $3000$ to obtain the first function $F_{0,1}$, and the number of training iterations to obtain the second function $F_{0,2}$ is set to $n_2 = 2  n_1.$

\begin{figure}
\begin{subfigure}[t]{0.5\textwidth}
\centering
    \includegraphics[width=\textwidth]{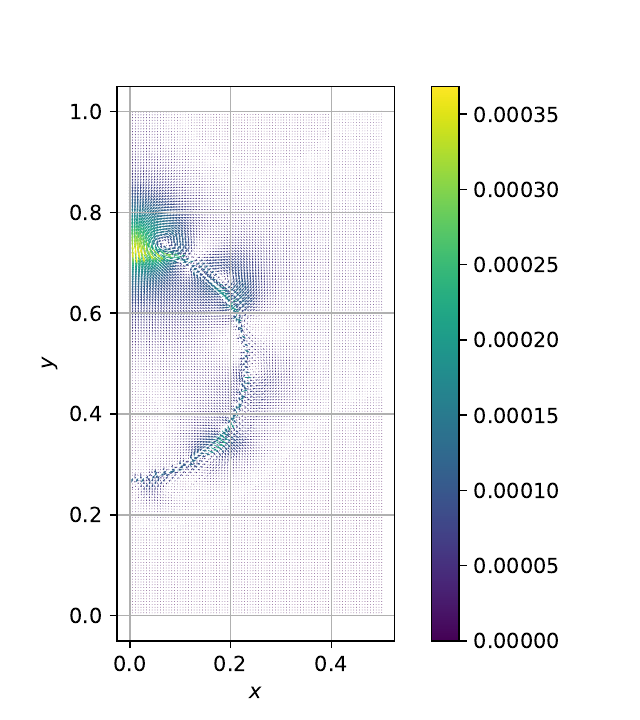}
\caption{}
\end{subfigure}
\begin{subfigure}[t]{0.5\textwidth}
    \includegraphics[width=\textwidth]{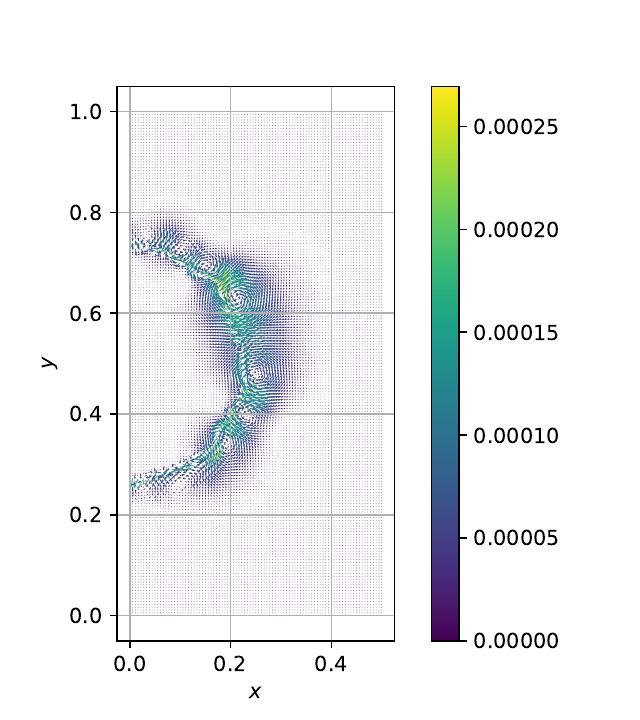}
\caption{}
\end{subfigure}
\begin{subfigure}[t]{0.5\textwidth}
    \includegraphics[width=\textwidth]{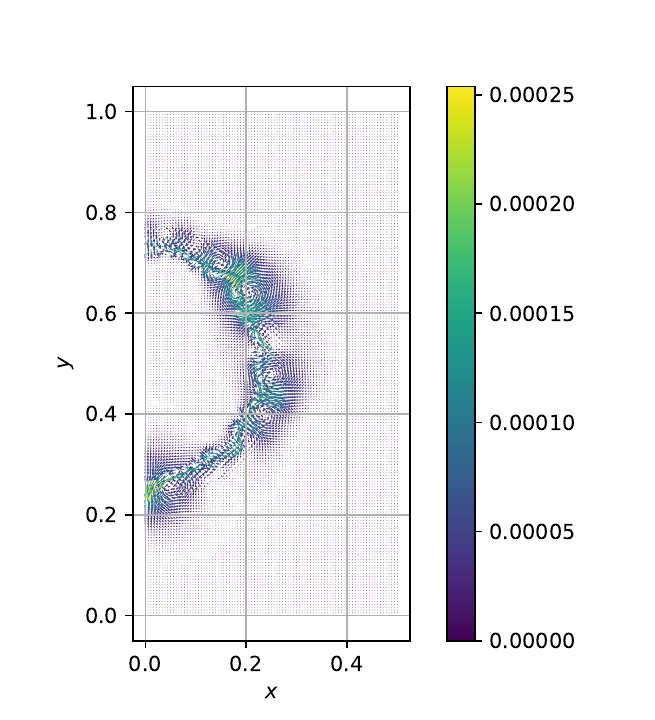}
\caption{}
\end{subfigure}
\begin{subfigure}[t]{0.5\textwidth}
    \includegraphics[width=\textwidth,height=\textwidth]{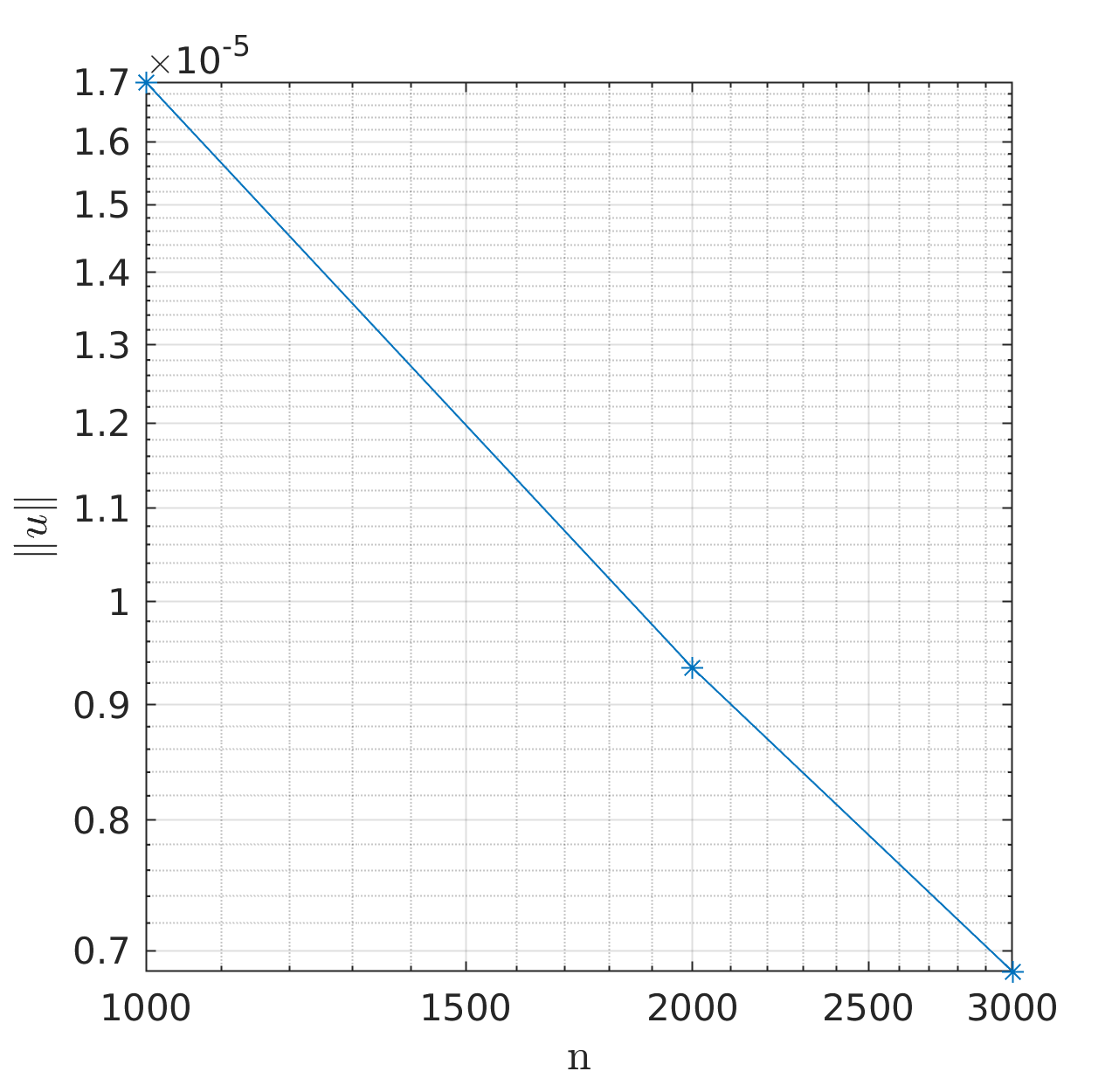}
\caption{}
\end{subfigure}
\caption{ The spurious velocities for the steady-state Laplace problem, obtained using three different 
number of training iterations, where the first function $F_{n,1}$ is obtained with (a) $n_1 = 1000$ training iterations, (b) $n_1 = 2000$ training iterations,(c) $n_1 = 3000$ training iterations. In (d) we show the $L2$ norm of the velocities for a different number of training iterations. The second function $F_{n,2}$ is obtained by training for $n_2 = 2n_1$ iterations.   }
\label{fig:V_level}
\end{figure}

The results shown in figure~\ref{fig:V_level} indicate that the approximation error of the PINN on the level-set function is an additional source for the creation of spurious currents, that can be diminished by increasing the number of training iterations.

\subsection{Rising bubble problem}

The next test is the 2D rising bubble problem, where we compare our results to the benchmarked results performed in \cite{BRbench}. The domain geometry is rectangular with a height of 2 and a width of 1. The initial conditions are $F_0 = r - 0.25$ for the level-set function, where $r$ is the distance to the point $[x,y]=[0.5,0.5]$, and the velocity field $\boldsymbol{u} = \boldsymbol{0}$. The physical parameters used in \cite{BRbench} can be found in Table \ref{table:1}. The non-dimensional numbers are calculated with the following characteristic values:
$\hat{L} = 1$;
$\hat{u} = 1;$
$\hat{\rho} = \rho_1$;
 $\hat{\mu} = \mu_1$.
With these characteristic values, we obtain that ${\rm Re} = 100$, ${\rm We}=1000/24.5$, and ${\rm Fr} = 1/\sqrt{0.98}$ for the first test case, and ${\rm Re} = 100$, ${\rm We}=1000/1.96$, and ${\rm Fr} = 1/\sqrt{0.98}$ for the second test case.

The results documented in~\cite{BRbench} are obtained using three different simulation codes. 
\begin{table}[ht!]
\centering
\begin{tabular}{ |c|c|c|c|c|c| } 
 \hline
 Test case & $\rho_1$ & $\rho_2$ & $\mu_1$ & $\mu_2$ & $\gamma$ \\ 
 \hline
 1 & 1000 & 100 & 10 & 1 & 24.5 \\ 
 2 & 1000 & 1 & 10 & 0.1 & 1.96 \\
 \hline
\end{tabular}
\caption{Physical parameters for the rising bubble test.}
\label{table:1}
\end{table}

\begin{figure}[ht!]
\centering
    \includegraphics[width=5.5cm]{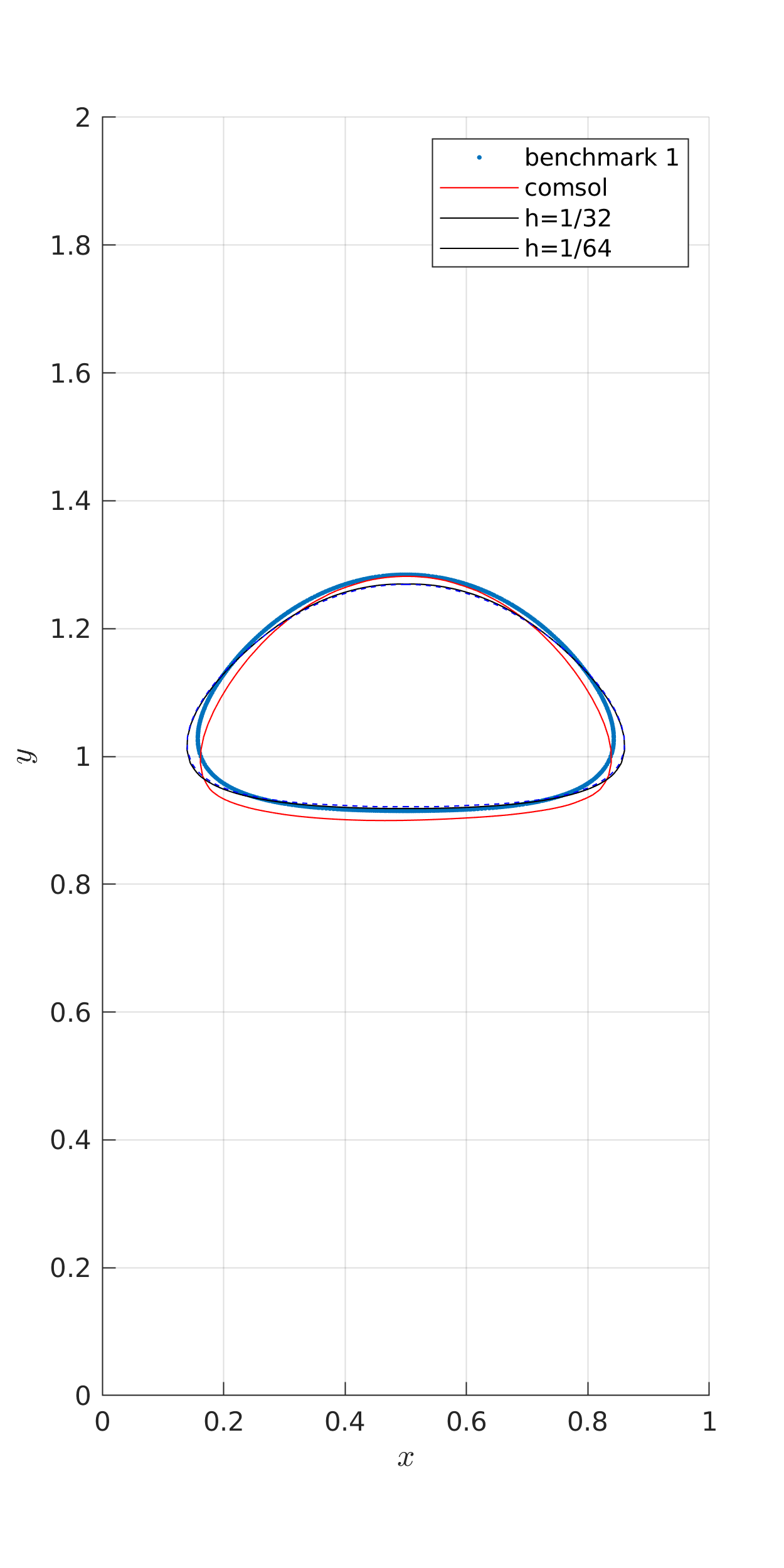}

    \caption{ The bubble shape obtained at $t=3$s for the first rising bubble test case problem, with the PINN-FEM coupling with elements size of $h=1/32$ and $h=1/64$, Comsol, and the results obtained in \cite{BRbench}.
    }
    \label{fig:bubble1}
\end{figure}

\begin{figure}[ht!]
\begin{subfigure}[t]{0.5\textwidth}
\centering
    \includegraphics[width=\textwidth]{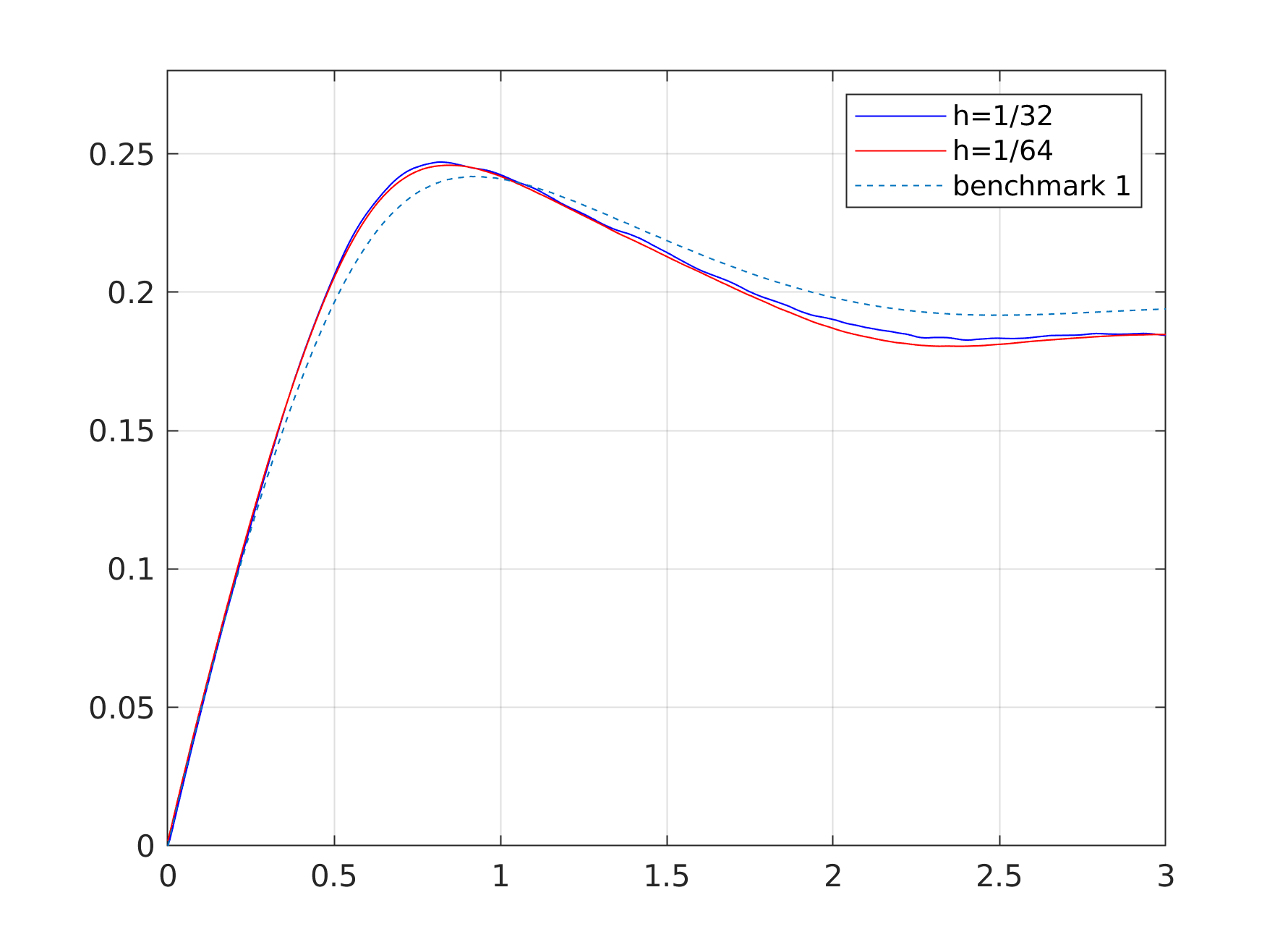}
\caption{}
\end{subfigure}
\begin{subfigure}[t]{0.5\textwidth}
    \includegraphics[width=\textwidth]{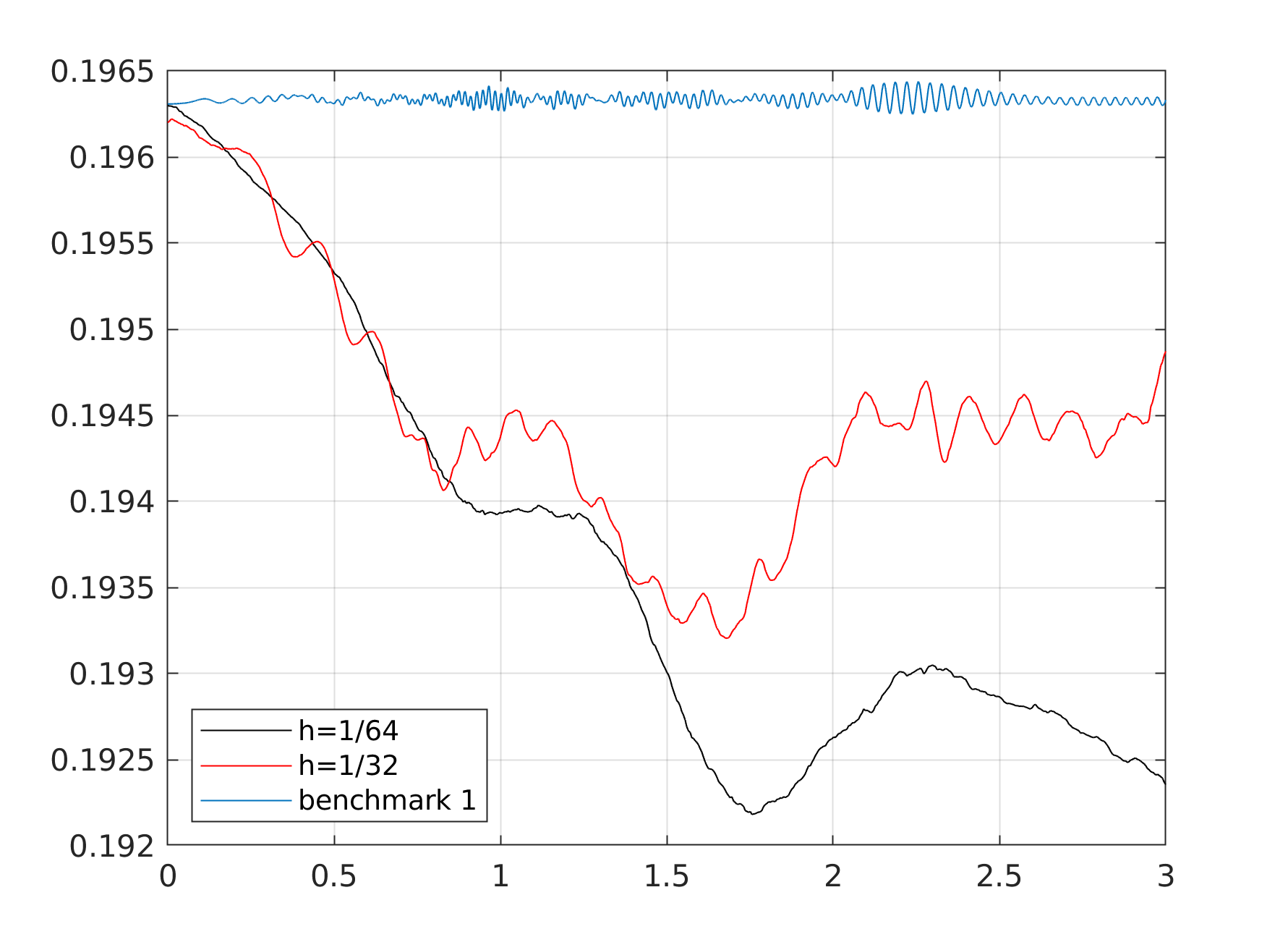}
\caption{}
\end{subfigure}
    \caption{ (a) The bubble rise velocity, and (b) the mass loss, for the first rising bubble test case problem, obtained with the PINN-FEM coupling with elements size of $h=1/32$ and $h=1/64$, and the results obtained in \cite{BRbench}. 
    }
    \label{fig:b1_r_m}
\end{figure}

For the first case, the bubble undergoes little deformation, and the codes tested by \cite{BRbench} gave the same bubble shape and rise velocity. The time step chosen for our simulation is ${\rm \Delta} t=0.005$. We reinitialize the level-set function at every time step, and we use 200 training iterations per time step to obtain the first function $F_{n,1}$ and 600 training iterations per time step to obtain the second function $F_{n,2}$, using the Adam optimizer. From the results presented in figures~\ref{fig:bubble1} and~\ref{fig:b1_r_m}, we notice good agreement with the results obtained by \cite{BRbench} and the Comsol solver. The mass loss is around $1.8 \% ,$ and the rise velocity of the bubble follows the same behavior as the velocities found in \cite{BRbench}.

\begin{figure}[ht!]
\begin{subfigure}[t]{0.6\textwidth}
\centering
    \includegraphics[width=5.5cm,height=11cm]{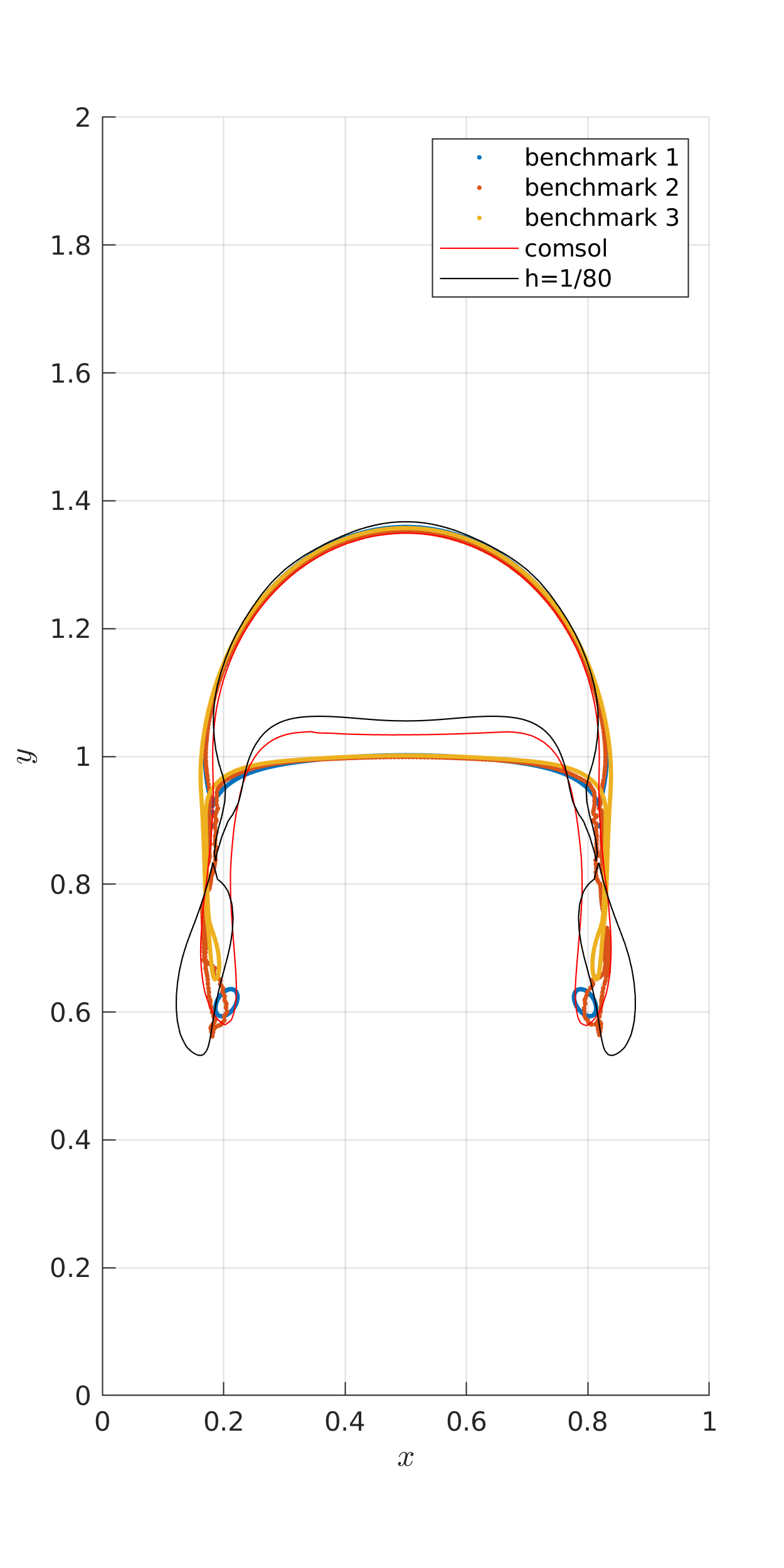}
\caption{}
\end{subfigure}
\begin{subfigure}[t]{0.35\textwidth}
    \includegraphics[width=5.5cm,height=11cm]{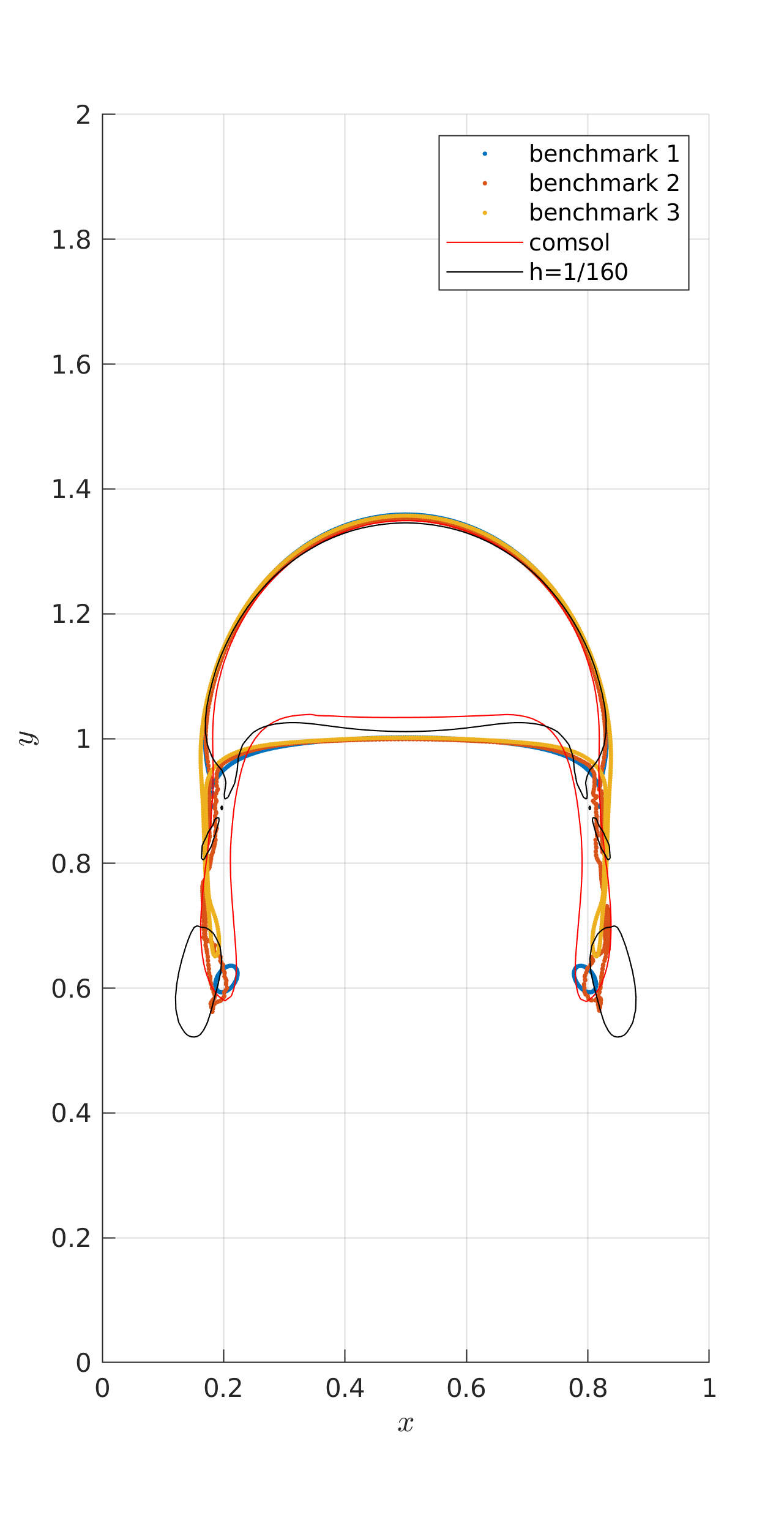}
\caption{}
\end{subfigure}
    \caption{The bubble shape obtained at $t=3$s for the second rising bubble test-case, obtained by the PINN-FEM coupling, (a) with elements size $h=1/80$, and (b) with $h=1/160$, Comsol, and the results obtained in \cite{BRbench}.
    }
    \label{fig:BR2shape}
\end{figure}

The second case is more challenging because the bubble will break, leaving behind thin filaments. The codes tested in \cite{BRbench} did not produce similar bubble shapes, particularly in the filaments and the broken-off bubbles. In our simulation, we use a time step of ${\rm \Delta} t=0.005$. We reinitialize the level set function every ten time steps, and we use 200 training iterations per time step to obtain the first function $F_{n,1}$, and 600 training iterations per time step to obtain the second function $F_{n,2}$, using the Adam optimizer.
The results are portrayed in figure \ref{fig:BR2shape}, where we can see that the proposed methodology yields comparable shapes to the benchmarked results.  

%\end{figure}

\begin{figure}[ht!]
\begin{subfigure}[t]{0.5\textwidth}
\centering
    \includegraphics[width=\textwidth]{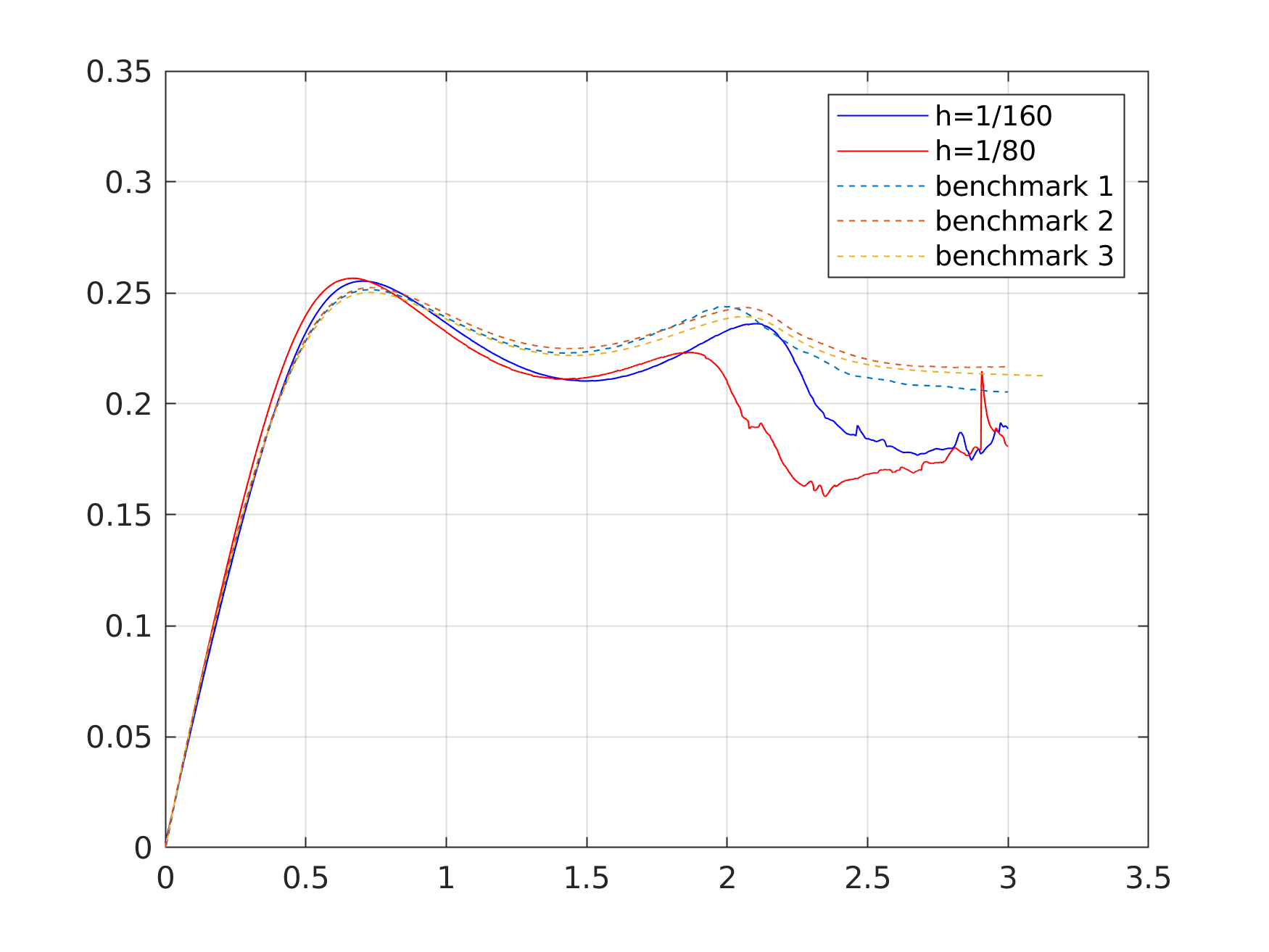}
\caption{}
\end{subfigure}
\begin{subfigure}[t]{0.5\textwidth}
    \includegraphics[width=\textwidth]{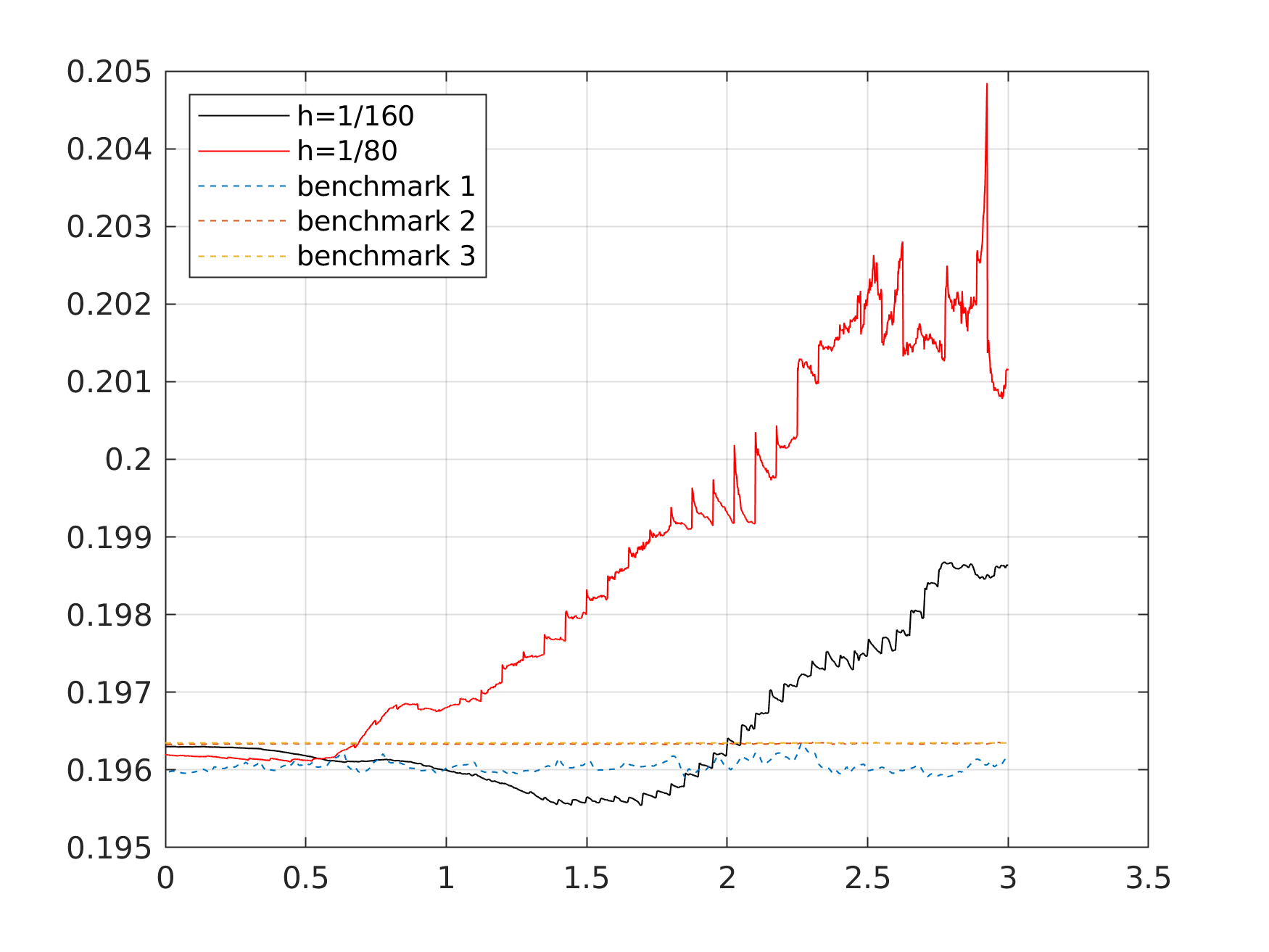}
\caption{}
\end{subfigure}
    \caption{ (a)The bubble rise velocity, (b) the mass loss, for the second rising bubble test case problem. The mass loss obtained is $2.5 \%$ with elements size $h=1/80$, and $1.2 \%$ with $h=1/160$.}
    \label{fig:BR2quant}
\end{figure}
% \subsection{Performance of the proposed methodology}
% The proposed methodology mixes PINN and the FEM. In the FEM, we obtain a more accurate approximation of the velocity field for smaller time steps, and the number of degrees of freedom determines the cost of computation per time step. 

% On the other hand, the accuracy of the approximated level-set function with PINN is not heavily affected by the time step size. The accuracy of the approximation and the computational cost are mainly governed by the number of training iterations per time step. 
% When using small time steps, the level-set function does not vary significantly between two consecutive time steps, this means that we can use a smaller number of training iterations per time step for smaller time steps.

\section{Conclusion}
We developed a coupled PINN/FEM method for the simulation of multi-fluid flows. The Navier-Stokes equations were discretized using the FEM, with the P2/P1 element combination for the velocity and pressure. The level-set function was approximated using a multi-level PINN, where the first PINN approximates the level-set function in the entire domain, and the second PINN refines the approximation near the interface. Three PINN formulations were proposed to handle the imposition of the initial conditions: weak imposition, strong imposition, and discrete-time steps. Numerical tests showed that the strong imposition of initial conditions gives the most accurate results and that the accuracy of the level-set function is mainly affected by the number of training iterations per time step. 
Three reinitialization strategies were proposed, that give a differentiable solution by automatic differentiation. The first two strategies find the distance from any point in the domain to the interface by constrained minimization. The third strategy reinitializes the level-set function using a PINN for the Eikonal equation, with a strong imposition of the interface position. Numerical tests showed that the reinitialization by a PINN yields more accurate results for the curvature calculations, and better mass-conservation properties. Finally, the proposed numerical strategy was tested against the rising bubble problem, benchmarked in \cite{BRbench}. It was found that the proposed strategy can produce results comparable to traditional methods.
\clearpage
\bibliographystyle{abbrv} 
\bibliography{pinn}

\end{document}